\documentclass[12pt]{article}
\usepackage{latexsym, amsfonts, amsmath}

\def\boxtimes{\Box\kern -6pt\times}

\newtheorem{thm}{Theorem}[subsection]
\newtheorem{lem}[thm]{Lemma}
\newtheorem{cor}[thm]{Corollary}
\newtheorem{prop}[thm]{Proposition}

\newtheorem{rem}[thm]{Remark}

\def\ovl{\overline} 
\def\h{\hbar}
\def\wt{\widetilde}
\def\O{{\mathcal O}} 
\def\reg{^{reg}}

\def\g{{\mathfrak g}}
\def\h{{\mathfrak h}}

\def\k{{\mathfrak k}}

\def\n{{\mathfrak n}}
\def\p{{\mathfrak p}}

\def\u{{\mathfrak u}}

\def\so{{\mathfrak so}}
\def\sp{\mathfrak{sp}}
\def\sl{\mathfrak{sl}}
\def\gl{\mathfrak{gl}}

\def\dimC{\dim}

\def\Hom{\mbox{\rm Hom}}
 
\def\Th{\Theta}
\def\bA{{\bf A}} \def\bB{{\bf B}} \def\bC{{\bf C}} \def\bS{{\bf S}}
\def\half{\frac{1}{2}}
\def\ep{\epsilon}
\def\sy{{\bf symbol}}
\def\la{\lambda}
\def\pdb#1{\frac{\partial\phantom{x}}{\partial#1}}
\def\C{{\mathbb C}}
\def\Z{{\mathbb Z}}
\def\R{{\mathbb R}}

\def\D{{\mathcal D}} \def\A{{\mathcal A}}
\def\bS{{\bf S}}
\def\wdg{\wedge}
\def\tw{\widetilde{w}} \def\tx{\widetilde{x}} \def\tf{\widetilde{f}} 
\def\bT{{\bf T}}
\def\U{{\mathcal U}}
\def\i{^{-1}}
 
\def\mypar{\vskip 0pc}
\def\heis{\g_{neg}}
\def\Oreg{\O^{reg}}
\def\End{\mbox{End}}  \def\Hom{\mbox{Hom}}
\def\ov{\ovl{v}}  \def\oV{\ovl{V}}

\def\og{\ovl{g}}  \def\of{\ovl{f}} \def\oD{\ovl{D}}
\def\RR{R_{-1}^{\phantom{i}4}}
\def\Cl{\mathrm{Cl}}

\def\Qed{$\Box$\medskip}
\def\ad{\mbox{ad}}  
\def\AND{\qquad\mbox{and}\qquad} 
\def\d{{\mathord{\rm d}}}
\def\ddto{\frac{\d\phantom{t}}{\d{t}}|_{t=0}}

\def\Vect{\mathfrak{Vect}}
\def\Omin{O_{\mathord{\rm  min}}}
 
\title{Exotic Differential Operators on Complex  Minimal Nilpotent Orbits
}

\author{Alexander Astashkevich\thanks{Department of Mathematics, 
           UC Davis, E-mail address: ast@math.ucdavis.edu} \, and
        Ranee Brylinski\thanks{Department of Mathematics, 
        Penn State University, E-mail address: rkb@math.psu.edu}  
\thanks{Research supported in part  by NSF  Grant No. DMS-9505055}}

\date{}
\begin{document} 
\maketitle

\begin{abstract} 
The orbit $\O$ of highest weight vectors in a complex simple Lie algebra $\g$
is the smallest non-zero adjoint orbit that is conical, i.e., stable under the 
Euler dilation action of $\C^*$ on $\g$. 
$\O$ is a smooth quasi-affine variety and its closure $\Cl(\O)$ is a singular 
quadratic cone in $\g$.   The varieties $\O$ and  $\Cl(\O)$ have the same algebra
of regular functions, i.e. $R(\O)=R(\Cl(\O))$, 
and hence the  same  algebra of algebraic differential operators, i.e.,
$\D(\O)=\D(\Cl(\O))$. 

The content and structure of $\D(\O)$ is unknown, except
in the case where $\g=\sp(2n,\C)$.  The  Euler action defines natural
 $\g$-invariant algebra gradings
$R(\O)=\oplus_{p\in\Z_+} R_p(\O)$ and $\D(\O)=\oplus_{p\in\Z}\D_p(\O)$
so that $D(R_q(\O))\subset R_{q+p}(\O)$ if  $D\in \D_p(\O)$.

We construct, for $\g$ classical and different from $\sp(2n,\C)$,  a  
subalgebra $\A$ inside
$\D(\O)$ which, like $R(\O)$, is  $\g$-stable,
graded and maximal commutative. Unlike $R(\O)$, which lives in non-negative
degrees,  our algebra $\A$ lives in non-positive degrees so that
$\A=\oplus_{p\in\Z_+}\A _{-p}$.  The intersection 
$\A\cap R(\O)$ is just the constants.   The  space 
$\A_{-1}$ is finite-dimensional, transforms as the
adjoint representation and generates $\A$ as an algebra.
Moreover we construct a $\g$-equivariant graded
algebra isomorphism  $R(\O)\to \A$, $f\mapsto D_f$. 
The differential operators in $\A_{-1}$  have order $4$. 

The  operators in $\A$ are ``exotic"  in that they lie outside the realm
of familiar differential operators.  Although their existence is predicted
(as a conjecture)  by the
quantization program for nilpotent orbits of the second author, our  actual
construction, by  quantizing ``exotic" symbols we obtained previously, is
quite subtle.  

We will use the  operators in $\A_{-1}$  to quantize $\O$ in a subsequent  paper
and to construct an algebraic star product on $R(\O)$. To pave the way for the
quantization, we show that the
formula $(f|g)=(D_{\og}\,f)(0)$ defines a  positive-definite
Hermitian inner product on $R(\O)$. 
\end{abstract}

\section{Introduction}
\label{sec:intro} 
 
The nature of the algebra $\D(X)$ of differential operators on a complex 
quasi-projective algebraic  variety $X$ is rather mysterious. There is only one 
rather  general result: if $X$ is both smooth and affine, then (i) $\D(X)$ is generated
by the functions and  vector fields on $X$, and (ii) every function on $T^*X$ which is
a homogeneous polynomial on the fibers of $T^*X\to X$ is the principal symbol of a
differential operator. As soon as $X$ becomes non-smooth or non-affine, the
generation quickly fails and the surjectivity of the symbol map is not well
understood.

Let $X$ be a homogeneous space of a reductive  complex algebraic 
group $K$.  Then $X$ is a quasi-projective complex algebraic manifold.
Differentiation of the $K$-action gives an infinitesimal
action of the Lie algebra $\k$ on $X$ by vector fields. This gives rise to
a complex algebra homomorphism  
\begin{equation}\label{eq:1.1}
\U(\k)\to\D(X)
\end{equation}
where $\U(\k)$ is the universal  enveloping algebra of $\k$.  

Suppose $X$ is  compact, i.e.,    projective. Then W. Borho and J.-L. Brylinski 
(\cite{bobr})  proved that the  homomorphism (\ref{eq:1.1})  is
surjective. They also computed its kernel which is   a primitive ideal in $\U(\k)$. 
Here $X$ is a (generalized) flag variety $K/Q$. These spaces $X$ occur as the
projectivized    orbits of highest weight vectors in finite-dimensional irreducible
representations of $K$.   We note that $K$ admits only finitely many compact
homogeneous spaces.  These orbits play a key role in  
algebraic geometry, symplectic geometry  and representation theory.

There are natural non-compact analogues to the flag varieties.
These homogeneous spaces   of $K$  arise   in the following way. 
Suppose that $K$ sits inside a larger reductive  complex  
algebraic  group $G$ as (the identity component of)
a spherical subgroup. Let $\g=\k\oplus\p$ be the corresponding  complex Cartan
decomposition where $\k=\mbox{Lie}\, K$ and $\g=\mbox{Lie}\, G$. Then the
$K$-orbits  in $\p$ are quasi-affine homogeneous spaces $X$ of $K$ (and so
non-compact except if $X$ is a point).  The choice $G=K\times K$ gives rise to all the
adjoint orbits of $K$.

Among all $K$-orbits  in $\p$, there are only finitely many orbits $Y$
which are stable under the Euler dilation action of $\C^*$ on $\p$.
These orbits are  exactly the ones consisting of nilpotent elements of $\g$;
accordingly they are called ``nilpotent orbits". There is a rich literature on the 
geometry of these orbits $Y$,  especially with respect to the quantization of 
the corresponding (see \cite{sek}) real nilpotent orbits. 
The most recent development here is the Kaehler structure on real nilpotent orbits 
arising from  the work of Kronheimer (\cite{kron}) and Vergne (\cite{ver}).
See also  \cite{brInstantons}, \cite{brGQmin}, \cite{brISC}.
 
The  Euler action on $Y$ provides $K$-invariant  algebra gradings  
$R(Y)=\oplus_{p\in\Z}R_p(Y)$ and  $\D(Y)=\oplus_{p\in\Z}\D_p(Y)$  
where $R_p(Y)$ is the subspace of Euler homogeneous functions of degree $p$  and   
\begin{equation}\label{eq:1.2}
\D_{p}(Y)=\{D\in\D(Y)\,|\,  D(R_q(Y))\subset R_{q+p}(Y)\,\,
\mbox{for all $q\in\Z$}\}
\end{equation}
We say that  differential operators in $\D_{p}(Y)$ are
{\it Euler homogeneous of   degree $p$}. The vector fields of the infinitesimal 
$\k$-action are  Euler homogeneous of   degree $0$ since the $K$-action
commutes with the Euler action.

$\D(Y)$ is not, in general,  generated by the global functions and vector fields 
on $Y$. We can see this already in the simplest case,  namely the case
where $\g=\sp(2n,\C)$,  $\k=\gl(n,\C)$ and $Y$ is  a minimal (non-zero)  orbit 
in $\p$. (So $Y$ is the orbit of highest weight vectors in $\p^+$ or $\p^-$.)
Then  $Y$ is isomorphic to the quotient  $(\C^n-\{0\})/\Z_2$ where  $\Z_2$ acts 
by $\pm 1$.  It is easy to see that $\D(Y)$ is the  even part of the Weyl algebra, i.e., 
\begin{equation}
\D(Y)=\textstyle
\C[x_1,\dots,x_n,\pdb{x_1},\dots,\pdb{x_n}]^{\Z_2}=
\C[x_ix_j, x_i\pdb{x_j},\frac{\partial^2\phantom{x}}{\partial x_i\partial x_j}]  
\end{equation}
Clearly the subalgebra generated by the functions and vector fields is
$\C[x_ix_j,x_i\pdb{x_j}]$. This subalgebra lives in
non-negative degrees, and so does not  contain the operators
$D_{ij}=\frac{\partial^2\phantom{x}}{\partial x_i\partial x_j}$. 
The operators $D_{ij}$ generate a maximal commutative subalgebra $\A$
of $\D(Y)$ which lives in non-positive degrees. 
 
In most cases beyond the example above,  it is a hard open problem to determine 
$\D(Y)$ and to construct any ``new operators" in it beyond the known operators, i.e.,
those  generated   by functions and by the vector fields of the $\g$-action.   See
Section \ref{sec:last} for a few  known examples where  highly
non-obvious  ``new " operators  in $\D(Y)$ have been constructed.

We expect $\D(Y)$ to contain   ``new"  differential  operators  for general
nilpotent orbits $Y$ (cf. \cite{brISC}).  These will be  ``exotic" in the sense
that they lie outside the realm  of familiar differential operators.

In this  paper, we  study the case where $Y$ is the complex minimal nilpotent 
orbit   $\O$ of a    simple complex Lie group $G$.   
(This is the case where we  embed  $K$ as a spherical subgroup of $K\times K$;
then we rename $K$ as $G$.)  See \cite{jos} and \cite{bkHam} for 
some results on the geometry and  quantization of $\O$.  

The space  $\O$ is just the orbit of highest weight vectors (or equivalently, of highest
root vectors) in $\g$. Accordingly $\O$ has a number of very nice properties.
The closure of $\O$ is $\Cl(\O)=\O\cup\{0\}$ and $\Cl(\O)$ is a quadratic cone
in $\g$, i.e., $\Cl(\O)$ is cut out (scheme-theoretically) by homogeneous quadratic
polynomial functions on $\g$.  All regular functions on $\O$ extend to
$\Cl(\O)$. So $\O$ and $\Cl(\O)$ have the same algebra  of regular functions, i.e.,
$R(\O)=R(\Cl(\O))$,  and hence the  same  algebra of algebraic differential operators, 
i.e.,  $\D(\O)=\D(\Cl(\O))$.  Notice that $\O$ ``misses" being a smooth affine variety 
by just one point, namely the origin.

The algebra $R(\O)$ is graded in non-negative degrees so
that $R(\O)=\oplus_{p\in\Z_+}R_p(\O)$. The subspace $R_1(\O)$ generates $R(\O)$
as an algebra and consists of  the functions $f_x$, $x\in\g$,  obtained  by restricting 
the linear functions on $\g$. So $f_x$ is defined by  $f_x(z)=(x,z)_\g$ where
$(\cdot,\cdot)_\g$ is the (normalized) Killing form of $\g$.

Differentiation of the $G$-action on $\O$ gives an infinitesimal vector field action 
\begin{equation}\label{eq:1-eta}
\g\to{\Vect}\,\O, \qquad   x\mapsto\eta^x  
\end{equation}
These vector fields $\eta^x$ are Euler homogeneous of degree $0$.
The subalgebra  of $\D(\O)$ generated by $R(\O)$ and by the  $\eta^x$
is thus Euler graded in non-negative degrees.  

We assume that $\g$ is classical in order to prove our results.
We exclude the case  $\g=\sp(2n,\C)$, $n\ge 1$, as then
$\O\simeq(\C^{2n}-\{0\})/\Z_2$ and so $\D(\O)$ is  just the even part of a Weyl
algebra.

We construct inside $\D(\O)$ a $G$-stable, graded, maximal commutative subalgebra 
$\A$ isomorphic to $R(\O)$. Unlike $R(\O)$, which lives in non-negative
degrees,  our algebra $\A$ lives in non-positive degrees so that
\begin{equation}\label{eq:A-intro} 
\A=\bigoplus_{p\in\Z_+}\A _{-p}
\end{equation}  
where $\A_{-p}=\A\cap\D_{-p}(\O)$.
The intersection $\A\cap R(\O)$ is just $\A_0=R_0(\O)=\C$.   The  space
$\A_{-1}$ is finite-dimensional, transforms as the
adjoint representation and generates $\A$ as an algebra.
We construct a $G$-equivariant graded algebra isomorphism  $R(\O)\to \A$, 
$f\mapsto  D_f$. 

The differential operators in $\A_{-1}$  have order $4$ and are ``exotic" in the sense
discussed above. We construct $\A_{-1}$  by quantizing a space 
$r(\g)\subset  R(T^*\O)$  of   symbols  which we obtained  in  \cite{abr2}.  
We next recall the context in which we found  $r(\g)$.

In  \cite{abr1}, we constructed   complex Lie  algebra embedding 
\begin{equation}\label{eq:tau-intro} 
\tau:\g\oplus\g\to\C(T^*\O)
\end{equation} 
where  $\C(T^*\O)$ is the complex Poisson  algebra of rational functions on $T^*\O$. 
We did this in a more general context which includes  the cases $Y=\O$.

The restriction of  $\tau$ to the diagonal in  $\g\oplus\g$ is the natural map
\begin{equation}\label{eq:Phi-intro} 
\Phi:\g\to R(T^*\O),\qquad x\mapsto\Phi^x
\end{equation} where $\Phi^x$ is the principal symbol of the 
vector field $\eta^x$.   I.e., $\tau(x,x)=\Phi^x$. On the anti-diagonal $\tau$ is given by
\begin{equation}\label{eq:b-intro} 
\tau(x,-x)=f_x-g_x
\end{equation} 
where   $g_x$, $x\in\g$,   are   rational functions on  $T^*\O$  defined on a
$G$-invariant (Zariski) open algebraic submanifold $(T^*\O)^{reg}$. 

The functions $g_x$  transform in the adjoint representation of $G$. Each function 
$g_x$ is homogeneous of degree $2$ on the fibers  of $T^*\O\to\O$. Also $g_x$ is 
homogeneous of degree $-1$ with respect to the canonical lift of the Euler scaling  
action on  $\O$.  The functions $g_x$, like the functions $f_x$,
Poisson commute; i.e., $\{f_x,f_y\}=\{g_x,g_y\}=0$ for all $x,y\in\g$.

In \cite{abr2}, we computed $g_x$ for a restricted class of orbits $Y$ which still
includes $Y=\O$. Let $\la\in R(T^*\O)$ be  the principal symbol of the Euler vector field
$E$ on  $\O$. We found that $(T^*\O)^{reg}$ is the  locus  in $T^*(\O)$ where
$\la$ is non-vanishing. We got the  formula:
\begin{equation}\label{eq:g_x} 
g_x=\frac{r_x}{\la^2}
\end{equation}
where
\begin{equation}\label{eq:r-intro} 
r:\g\to R(T^*\O),\qquad x\mapsto r_x
\end{equation}
is an explicit  complex linear $\g$-equivariant map.  

Each  function $r_x$ is  a homogeneous degree $4$ polynomial  on the fibers of
$T^*(\O)\to\O$. Also $r_x$ is homogeneous of degree $-1$
with respect to the  canonical lift of the Euler   action. The functions $r_x$ again
Poisson commute, i.e.,  $\{r_x,r_y\}=0$.

We say that a function $\phi\in R(T^*X)$ which is homogeneous of degree $d$
on the fibers of $T^*X\to X$ (and hence a polynomial function on the fibers)
is an  \emph{order $d$ symbol}.  The order $d$ symbol map sends 
$\D^d(X)$ into $R^d(T^*X)$ where $\D^d(X)$ is the  space  of differential operators of
order at most $d$ and $R^d(T^*X)$ is the  space  of order $d$ symbols. 
We say that $\phi\in R^d(T^*X)$ \emph{quantizes} into $D$ if $D\in\D^d(X)$
has order $d$  and the principal symbol of $D$ is  $\phi$.
Unless $X$ is smooth and affine, there is no guarantee that a given symbol will quantize.
 
In this paper we quantize the symbols $r_x$ into order $4$ differential
operators $D_x$ on $\O$ in a manner  equivariant with respect to both the  $G$-action
and the Euler $\C^*$-action. We show that this equivariant quantization is unique.
In our next paper \cite{abr3}, we use these same  operators $D_x$ to quantize
$\O$ by quantizing the map (\ref{eq:tau-intro}). We obtain a star product on
$R(\O)$ given by ``pseudo-differential" operators.

We construct the operators $D_x$ by manufacturing a single operator
$D_0=D_{x_0}\in\D^4_{-1}(\O)$ where $x_0\in\g$ is a lowest weight vector.
Here  $\D^d_{p}(\O)=\D^d(\O)\cap\D_{p}(\O)$.
We build $D_0$ so that its principal symbol is $r_0=r_{x_0}$ and $D_0$ is a lowest
weight vector for a copy of $\g$ in $\D^4_{-1}(\O)$.
Then all other operators $D_x$ are simply obtained by taking (iterated) commutators
of $D_0$ with the vector fields $\eta^y$.

To construct $D_0$, we use our explicit formula for $r_0$ from \cite{abr2} to
quantize $r_0$ into a differential operator on a (Zariski) open set $\Oreg$ in $\O$.
Our quantization method is a modification of Weyl quantization by symmetrization.
We are forced to modify the  quantization of a factor $\la^2$ appearing in a term
of $r_0$ in order that the quantization of $r_0$ is  $\g$-finite and hence extends to
an operator on all of $\O$.  See Subsection \ref{ssec:strat} for a detailed discussion of
our strategy.  

In essence, we are adding lower order correction terms to the
true Weyl quantization of $r_0$ in order to obtain $D_0$.
The correction is uniquely determined, but its  nature is mysterious to us.

This paper can  largely be read independently of \cite{abr1} and \cite{abr2}, as the
symbols  in fact only motivate the construction of the differential operators.  Once we
figure out the correct formula for $D_0$, we give a self-contained proof  that 
$D_0\in\D^4_{-1}(\O)$.  The more abstract and general results we prove in Section
\ref{sec:diffop} for  differential operators on cones of highest weight vectors then
give in particular the main properties of our operators $D_x$:
(i) the operators $D_x$  commute, (ii) the operators $D_x$ 
generate a maximal commutative subalgebra of $\D(\O)$, and (iii)
$\ovl{f_x}$ and $D_x$ are adjoint operators on $R(\O)$
with respect to a (unique) positive definite Hermitian inner  product $(\cdot|\cdot)$  on
$R(\O)$ such that $(1|1)=1$.

A basic question in quantization theory is whether one can set up a preferred
quantization converting functions on $T^*X$ which are polynomial on the fibers of
$T^*X\to X$ into differential operators. This certainly involves  quantizing 
symbols into operators which include ``lower order" terms. This has been studied in a 
different context in  \cite{lo} -- see  Section \ref{sec:last}.

In Section \ref{sec:last} we compare our  constructions  with the results of 
Levasseur, Smith  and Stafford  on the Joseph ideal and on differential operators  
on classical  rings of invariants, and with the constructions in \cite{bkDiff}.

\section{Differential operators on conical orbits}
\label{sec:diffop}

\subsection{Differential operators on $G$-varieties}

Let $X$ be a complex algebraic quasi-projective variety.  Let $\D_X$ be the  
sheaf of algebraic differential operators on $X$. 
Then $\D(X)=\Gamma(X,\D_X)$ is the algebra of (global) differential 
operators  on    $X$.  Let $\D^d(X)\subset\D(X)$ be the subspace of differential 
operators of order at most $d$.  See e.g.  \cite[Appendix]{bkDiff} for  some  of the
basic definitions and facts on differential operators.

Let $R(X)$ denote the algebra of regular functions on $X$. For any 
commutative  complex algebra $A$, let ${\mathfrak D}(A)$ be the algebra  
of differential operators of $A$. We have a natural algebra homomorphism
$\D(X)\to {\mathfrak D}(R(X))$ and this is an isomorphism if $X$ is
quasi-affine. We will be dealing primarily  with differential operators on quasi-affine 
varieties.

Throughout the paper, we are working  with 
complex quasi-projective algebraic   varieties. 
The  functions,  vector fields, and differential operators we deal with are all regular, 
or algebraic, in the sense of  algebraic geometry.   

Our main emphasis is on working with smooth varieties, which are then
complex algebraic manifolds and admit cotangent bundles.    We often
use the term ``algebraic holomorphic" in place of ``regular" to emphasize that the
algebraic structure is an overlay to the underlying holomorphic structure.

Now suppose $X$ is smooth.  Let 
$G$ be a reductive  complex algebraic group. Assume $G$  acts on
$X$ by an  algebraic holomorphic group  action. Then there are induced representations 
of $G$   on $R(X)$  and  on $\D(X)$ by (complex) algebra automorphisms. Both
$G$-representations are  locally  finite and  completely reducible.

We use the term {\it $G$-linear map} to mean a $G$-equivariant complex linear map of 
vector spaces.  A   vector space $V$ is {\it $G$-irreducible} if $V$ is an irreducible
representation   of $G$; often $V$ will be a $G$-stable space inside a larger
representation.  If $V$ is $G$-irreducible and $W$ is some $G$-representation, then we  
say that a non-zero $G$-linear map $V\to W$ gives a ``copy of $V$ inside $W$". 
 
The differential of
the $G$-action on $X$ is the infinitesimal algebraic holomorphic vector field action  
\begin{equation}\label{eq:2-eta}
\g\to{\Vect}\,X, \qquad   x\mapsto\eta^x  
\end{equation}
where $\g=\mbox{Lie}(G)$ and $\eta^x_z=\ddto(\exp -tx)\cdot z$. These vector 
fields give  a complex Lie algebra homomorphism 
\begin{equation}\label{eq:2-etacomm}
\g\to\End \,\D(X), \qquad   x\mapsto\ad\,\eta^x  
\end{equation}
where $\ad\,\eta^x(D)=[\eta^x,D]$. This
representation of $\g$ on $\D(X)$ extends  canonically to a representation  of $\g$ 
on  $\D(X')$  where $X'$ is any Zariski open set in $X$. 

\begin{prop}
\label{prop:diffextend}
Suppose $G$ acts transitively on $X$.
Let  $D$ be a differential operator on a  Zariski open   set in $X$. Then $D$ extends 
to a  differential  operator  on $X$ if and only if $D$ is $\g$-finite. 
\end{prop}
{\bf Proof}: We have a filtration of $\D_X$ by the sheaves $\D_X^k$ of 
differential operators on $X$ of  order at most $k$. For each $k$,  $\D^k_X$ 
is the   sheaf of sections of some complex vector bundle on $X$ of finite rank.
So  the result follows by the same general lemma we used in \cite{abr2} to
prove the corresponding statement for $R(T^*X)$.
\Qed

We say that $X$ is a $G$-cone if there is a (non-trivial) algebraic holomorphic group 
action of 
$\C^*$ on $X$ commuting with the $G$-action.  The infinitesimal generator of the 
$\C^*$-action is the algebraic holomorphic {\it Euler vector  field}  $E$ on $X$. 
The Euler action defines $G$-invariant   algebra gradings  
\begin{equation}\label{eq:2-grad}
R(X)=\bigoplus_{p\in\Z}R_p(X)\AND \D(X)=\bigoplus_{p\in\Z}\D_p(X)
\end{equation}
where
\begin{eqnarray}
\label{eq:DpX} 
R_p(X)&=&\{f\in R(X)\,|\, E(f)=pf\}\nonumber\\[4pt]
\D_p(X)&=&\{D\in\D(X)\,|\, [E,D]=pD\}\nonumber
\end{eqnarray} 
Then, if  $X$ is quasi-affine,
\begin{equation}\label{eq:2-D-grad}
\D_p(X)=\{D\in\D(X)\,|\, D(R_q(X))\subset R_{p+q}(X) \mbox{ for all $q\in\Z$}\}
\end{equation}

The Euler grading of $\D(X)$ is compatible with the order filtration so that
we get a grading of the space of differential operators of order at most $d$
\begin{equation}\label{eq:D^d_p} 
\D^d=\bigoplus_{p\in\Z} \D^d_p
\end{equation}
where $\D^d_p=\D^d\cap \D_p$.

\subsection{Differential operators on the cone of highest weight vectors.}
\label{ss:2.2}

In this subsection, we assume that $X$ is the $G$-orbit of 
highest weight vectors in  a (non-trivial) finite-dimensional irreducible complex
representation $V$ of  $G$.  
A (non-zero) vector $v$ in $V$ is  called a {\it highest weight  vector} 
if $v$ is semi-invariant under some Borel subgroup $B$ of $G$.
(In Cartan-Weyl highest weight theory, one usually fixes a choice of $B$, and then
there is, up to scaling, exactly one $B$-semi-invariant vector. The conjugacy of Borel
subgroups implies that the set of all highest weight vectors is an orbit.) 

$X$ is a locally closed  subvariety of $V$ and   is stable under 
the scaling   action of  $\C^*$ on $V$. So $X$ is a $G$-cone with Euler 
vector  field $E$. 
In fact   the quotient   $X/\C^*$ is  the unique closed $G$-orbit in the projective
space ${\mathbb  P}(V)$. The quotient   $X/\C^*$is  often called the projectivization
${\mathbb P}(X)$ of $X$.    The closure  $\Cl(X)$ of   $X$ in 
$V$ is equal to $X\cup\{0\}$ and $X$ is the unique non-zero minimal conical
$G$-orbit in $V$.

The orbit $X$ has several very nice properties (\cite{pv}):
\begin{description}
\item[(i)] the closure $\Cl(X)$ is normal.
\item[(ii)] the graded $G$-equivariant algebra homomorphism $R(\Cl(X))\to R(X)$
given by restriction of functions  is an isomorphism.
\item[(iii)] the algebra $R(X)$ is multiplicity-free as a
$G$-representation and  $R_p(X)$ is $G$-isomorphic to the $p$th Cartan power of 
$V^*$.
\end{description}
These properties of $X$ can be derived using the Borel-Weil theorem on  
${\mathbb  P}(X)$.  

If a  $G$-irreducible vector space $W$ carries the  $G$-representation $V_{\mu}$, 
then the {\it $p$th Cartan power of $W$} is the (unique) $G$-irreducible subspace 
$W^{\boxtimes p}$ of the  $pth$ symmetric  power $S^p(W)$ which carries
$V_{p\mu}$.     Here $V_{\nu}$ denotes the finite-dimensional irreducible
$G$-representation of highest weight $\nu$. 
Geometrically, if $W$ is the space of  holomorphic sections of a $G$-homogeneous
complex line bundle $\cal L$ over ${\mathbb P}(X)$, then $W^{\boxtimes p}$ identifies 
with the space of holomorphic sections of  ${\cal L}^{\otimes p}$.

The algebra $R(V)$ identifies with the symmetric algebra $S(V^*)$.
Properties (ii) and (iii)  say  that  
the natural graded algebra homomorphism 
\begin{equation}\label{eq:S(V^*)} 
\zeta:S(V^*)\to R(X)
\end{equation}
defined by restriction of polynomial functions from $V$ to $X$ is surjective and
the  $p$th graded map $\zeta_p:S^p(V^*)\to R_p(X)$ induces a $G$-linear
isomorphism $\zeta_p:(V^*)^{\boxtimes p}\to R_p(X)$.

Consequently,  $R_p(X)=0$ for $p<0$ and 
\begin{equation}\label{eq:R_0}
R_0(X)=\C
\end{equation}
Also  $R(X)$ is generated by $R_1(X)$.   In degree $1$, $\zeta$ gives the $G$-linear
isomorphism 
\begin{equation}
\label{eq:fv}
V^*\to R_1(X),\qquad \alpha\mapsto f_\alpha
\end{equation}
where $f_\alpha(w)=\alpha(w)$.   Also property (ii) implies that
\begin{equation}
\label{eq:DCLX}
\D(X)=\D(\Cl(X))
\end{equation}

We have a natural graded algebra inclusion of $R(X)$ into $\D(X)$ as the space of
order zero differential operators; so  $R_p(X)\subset\D ^0_p(X)$.
Moreover $R(X)$ is a  maximal commutative  subalgebra of $\D(X)$
living  in non-negative degrees.   A natural
question is to find other maximal commutative subalgebras, and to see what degrees
they live in.

To figure out what we might expect to find, let us consider the simplest example.
This occurs when $V=\C^N$ is the standard representation of 
$G=GL(N,\C)$ and so  $X=\C^N-\{0\}$. Then $R(X)=R(\C^N)$ and so we get 
the  Weyl algebra 
\begin{equation}
\D(X)=\D(\C^N) =
\textstyle\C[v_1,\dots,v_n,\pdb{v_1},\dots,\pdb{v_n}]
\end{equation}
In addition to  $R(X)=\C[v_1,\dots,v_n]$, we have a
second ``obvious" maximal commutative subalgebra of $\D(X)$, namely  the algebra
\begin{equation}\label{eq:A-cc} 
\A=\C\textstyle[\pdb{v_1},\dots,\pdb{v_n}] 
\end{equation}
of constant coefficient differential
operators.   $\A$ lives in non-positive degrees.  The assignment
$v_i\mapsto\pdb{v_i}$  defines a graded algebra isomorphism $R(X)\to\A$.
$\A$ is generated by the space $\A_{-1}$ of constant coefficient vector fields
and $\A_{-p}\simeq V^{\boxtimes p}$ as $G$-representations.

We see no obvious way to generalize the construction of 
$L=\A_{-1}$ in the example above to the case where $V$ is arbitrary. 
However, if one has such a space $L$ then  we can prove
the following. This abstract  result  gives us the commutativity of the 
concrete exotic differential  operators we  construct later in Section 
\ref{sec:construct}.

\begin{prop} 
\label{prop:Vcomm}
Let $X\subset V$ and $X^*\subset V^*$ be the $G$-orbits of highest weight 
vectors. Suppose we have a  finite-dimensional $G$-irreducible complex subspace 
$L\subset\D_{-1}(X)$  which carries the $G$-representation $V$.  Then  the
differential operators $D\in L$ all commute
and hence generate a  commutative  subalgebra $\A\subset\D(X)$.
Consequently  $\A=\oplus_{p\in\Z_+}\A_{-p}$ is graded in non-positive degrees 
with  $\A_0=\C$ and $\A_{-1}=L$.

$\A$ is  isomorphic to $R(X^*)$. 
In fact, if we fix a $G$-linear isomorphism $V\to L$, $v\mapsto D_v$,  then 
the  map $R_1(X^*)\to\A_{-1}$, $f_v\mapsto D_v$, extends uniquely  to a
graded   $G$-linear  algebra  isomorphism   $R(X^*)\to\A$, $f\mapsto D_f$. 
\end{prop} 
{\bf Proof:}
Let us consider the complex $G$-linear   map 
$\pi:\,\wedge^2V\to \End\, R(X)$
given by $\pi(u\wdg v)=[D_u,D_v]$.  The image of $\pi$ lies in the space
$\End_{[-2]} R(X)$ of endomorphisms which are homogeneous of degree 
$-2$  with respect to  the   Euler action.  So $\pi$ is just a collection
of complex  $G$-linear maps
\begin{equation}
\label{eq:wedgemap}
\pi_p:\wedge^2V\to \Hom(R_p(X), R_{p-2}(X)) 
\end{equation}
Thus $\pi_p$ is an intertwining operator between  two $G$-representations, 
namely  $\wedge^2 V$ and  
$\Hom((V^*)^{\boxtimes p}, (V^*)^{\boxtimes(p-2)})\simeq 
\Hom(V^{\boxtimes(p-2)}, V^{\boxtimes p})$.

But it follows from  highest weight theory  that
\begin{equation}
\label{eq:Hom}
\Hom_{G}\left(\wedge^2 V,\Hom(V^{\boxtimes(p-2)},V^{\boxtimes p})
\right)=0.
\end{equation}
Indeed, let $\nu$ be the highest weight of $V$ so that 
$V^{\boxtimes k}\simeq V_{k\nu}$.
Every irreducible $G$-representation occurring in  $\wedge^2 V$ has
highest weight $\la<2\nu$ by Cartan product theory. On the other hand, 
suppose 
$\mu$ is the  highest weight of some  irreducible representation occurring in 
$\Hom(V_{(p-2)\nu},V_{p\nu})$. The theory of Parasarathy, Rao and 
Varadarajan for tensor product decompositions  gives here that 
$\mu\ge -(p-2)\nu+p\nu=2\nu$.  Thus $\la=\mu$ is impossible. This proves 
(\ref{eq:Hom}).

Therefore $\pi_p=0$ for all $p$ and so $\pi=0$. This proves the 
commutativity  of the operators  $D_v$.  The chosen map $V\to L$ defines a 
surjective $G$-linear graded
algebra homomorphism $S(V)\to\A$.   The image of $S^p(V)$ is the
graded component $\A_{-p}$.   Clearly $\A_{-p}\neq 0$ for $p\ge 0$.

We claim that the two graded  algebra
homomorphisms  $S(V)\to\A$ and $S(V)\to R(X^*)$ have the  same kernel.  
Indeed, suppose a $G$-irreducible subspace $W\subset \A_{-p}$ carries the 
representation  with  highest  weight  $\mu$. Then arguing as above, we find 
that  $\mu\le p\nu$ since $W$ occurs in $S^p(V)$ and  also   $\mu\ge p\nu$ 
because   each  operator $A\in\A_{-p}$ is a collection  of maps
$A_{j}:R_j(X)\to R_{j-p}(X)$.  This forces  $\mu=p\nu$.  But the representation 
$V_{p\nu}$ occurs exactly once in $S^pV$, as the Cartan piece $V^{\boxtimes p}$.  
We conclude that the space $V^{\boxtimes p}$  maps isomorphically to $\A_{-p}$.   
The claim now follows.\Qed

\begin{rem}\rm (i) Proposition \ref{prop:Vcomm} was established  in   
\cite[Theorem 3.10]{bkDiff} for a restricted class of cases $(\g,V)$.  We 
have essentially just rewritten the proof  given there in our more general 
situation.
\\ (ii)  Let $I\subset S(V)$ be  the graded ideal of polynomial functions vanishing on 
$X^*$. Then  we have a direct sum $S^p(V)=V^{\boxtimes p}\oplus I_p$.
An  important  property of 
$X^*$ is that $I$ is  generated by its degree two piece $I_2$. This is a result of 
Kostant; see \cite{gar} for a  write-up. This gives another way to prove the second  
paragraph  of Proposition  \ref{prop:Vcomm}.
\\(iii)  Suppose $L$ is as in   Proposition \ref{prop:Vcomm} but $L$ carries the 
representation  $V^*$ instead of $V$. Then we can show by a similar argument that  
$V\simeq V^*$ and so the conclusions of the Proposition still follow.
\end{rem}
 
Proposition  \ref{prop:Vcomm} leads us to pose the question:
in what generality does the space $L$    exist?  In \cite{bkDiff},  $L$ was constructed 
for    a restricted set of cases of highest weight orbits (associated to complex  Hermitian
symmetric pairs). The main result of this paper is to construct $L$ for the
minimal nilpotent orbit $\O$ in a classical complex simple Lie algebra  different from
$\sp(2n,\C)$. 

Our final  result of this subsection says, under a mild hypothesis, that the commutative 
subalgebra
$\A$ in Proposition \ref{prop:Vcomm} is maximal and the operators $D_f$  define a
non-degenerate inner product on $R(X)$. We will see later that our operators on $\O$ 
satisfy the hypothesis and the resulting inner product is positive definite.

Let $U\subset G$ be a maximal compact subgroup  of $G$. Then $V$ admits 
a  $U$-invariant positive-definite Hermitian inner product $(\cdot|\cdot)$  (unique 
up  to a positive real scalar).   Let $\oV$ denote the complex conjugate  vector space  
to $V$.  Then the complex conjugation map $V\to\oV$, $v\mapsto\ov$, is 
$\C$-anti-linear. Each  vector $\ov\in\oV$ defines a $\C$-linear
functional $b_{\ov}$ on $V$ by $b_{\ov}(u)=(u|v)$. This gives a  
$U$-equivariant identification  $\oV\to V^*$,
$\ov\mapsto b_{\ov}$, of complex vector spaces;  we will use this identification
freely from  now on.  This induces a graded $\C$-algebra  identification 
$\ovl{R(X)}=R(X^*)$. Then we get a  graded
$\C$-anti-linear algebra isomorphism  $R(X)\to R(X^*)$,   $g\mapsto\og$,
defined in degree $1$ by $f_{\ov}\mapsto\ovl{f_{\ov}}=f_v$.

\begin{prop}\label{prop:max}
Suppose we are in the situation of Proposition \ref{prop:Vcomm}. Let us assume
that for each positive integer $p$ there exists an operator    
$D\in L$  such  that  $D$ is non-zero on $R_p(X)$. Then 
\item{\rm(i)} The algebra $\A$, like $R(X)$, is a maximal commutative 
subalgebra  of  $\D(X)$.
\item{\rm(ii)}   $R(X)$ admits a unique  $U$-invariant non-degenerate
Hermitian inner product $(\cdot|\cdot)$, such that $(1|1)=1$  and, for all $v\in V$,  
the operators
$f_{\ov}$  and $D_v$ are adjoint on $R(X)$. Then for any $f,g\in R(X)$ we have
\begin{equation}\label{eq:HIP}
(f|g)=\mbox{constant term of } D_{\og}(f)
\end{equation}
\end{prop}
Each $h\in R(X)$ can be uniquely written as a finite sum $h=\sum_{p\in\Z_+}h_p$
where $h_p\in R_p(X)$. Then, in view of (\ref{eq:R_0}),  we call $h_0$ the          
{\it constant term} of  $h$.

{\bf Proof:}
We see that (ii) implies (i) since $R(X)$ is a maximal commutative subalgebra of  
$\D(X)$. To prove (ii), we need to construct $B=(\cdot|\cdot)$.

We start  by taking any  $U$-invariant
positive-definite  Hermitian inner product $Q$ on $R(X)$ such that 
$Q(1,1)=1$. Clearly  $Q$ exists and $R_p(X)$ and $R_q(X)$ are
$Q$-orthogonal if $p\neq q$ (since they carry different irreducible 
$U$-representations). Thus $Q$ is given by a family, indexed by $p$, of inner
products $Q_p$ on $R_p(X)$. Similarly $B$ must be given by a family $B_p$.

To start off, we  put  $B_0=Q_0$.  Now we proceed by induction and define
$B_{p+1}$ by the relation
\begin{equation}
B_{p+1}(h, f_{\ov}\, g)=B_p(D_v(h), g)
\end{equation}
where $g\in R_p(X)$ and  $h\in R_{p+1}(X)$.
This relation is exactly the condition that multiplication by $f_{\ov}$ is adjoint to
$D_v$.
 
We need to check that $B_{p+1}$ is  well-defined. Clearly the functions 
$f_{\ov}g$ span $R_{p+1}(X)$. Also  there exists a complex scalar 
$c_{p+1}$  such that
\begin{equation}
c_{p+1}Q_{p+1}(h, f_{\ov}\,g)=B_p(D_v(h),g)
\end{equation}
To see this we observe that  the two assignments
$v\otimes\ovl{g}\otimes h\mapsto Q_{p+1}(h, f_{\ov}\,g)$
and
$v\otimes\ovl{g}\otimes h\mapsto B_p(D_v(h),g)$ both define $U$-equivariant 
complex linear maps
\begin{equation}
R_1(X)\otimes \ovl{R_{p}(X)}\otimes R_{p+1}(X)\to\C
\end{equation}
But the space of such maps is $1$-dimensional -- this follows 
using highest weight  theory as in the
proof  of Proposition \ref{prop:Vcomm}.
So $B_{p+1}=c_{p+1}Q_{p+1}$. Thus $B_{p+1}$ is 
well-defined. (This was the same proof given in \cite[Theorem 4.5]{bkDiff}.)

Our hypothesis that  some $D_v$ is non-zero on $R_{p+1}(X)$
ensures  that  $c_{p+1}\neq 0$.   Hence $B_{p+1}$ is non-degenerate. 
Finally (\ref{eq:HIP}) follows as   $(f|g)=(D_{\og}(f)|1)$.   \Qed

\subsection{Complex minimal nilpotent orbits}
\label{ssec:2.3}

>From now on, we assume that $G$ is  a complex simple Lie group. Then   
the  Lie algebra  $\g=\mbox{Lie } G$ is a complex simple Lie algebra and so 
the  adjoint representation of $G$ on $\g$ is irreducible. We put $V=\g$
and  identify $V\simeq V^*$ in an $G$-linear way using the Killing form.

The $G$-orbit $X$ of highest weight vectors in $\g$ is now the 
minimal (non-zero) nilpotent orbit  $\O$ in $\g$.
Indeed,  we can fix a
triangular decomposition $\g=\n^+\oplus\h\oplus\n^-$ where $\h$ is a
complex Cartan  subalgebra and $\n^+$ and $\n^-$ are the spaces spanned
by, respectively,  the  positive and negative root vectors. Let $\psi\in\h^*$
be the highest (positive) root and let 
$H_\psi,Z_\psi,Z_{-\psi}$ be the corresponding root triple so that 
$H_\psi\in\h$  is semisimple and  $Z_ {\pm\psi}\in\n^{\pm}$ are nilpotent.
Then $Z_\psi$ is  a highest root vector in $\g$ and also a highest weight 
vector in  $\g$. So  $\g=V_{\psi}$ and 
$X$ is equal to the $G$-orbit $\O$ of  $Z_{\psi}$.

We  choose a real Cartan involution $\sigma$ of $\g$  such that $\h$ is 
$\sigma$-stable.  Then the $\sigma$-fixed space in $\g$ is a maximal compact Lie
subalgebra $\u$. Then  $\g=\u+i\u$ and  $\sigma$ is $U$-invariant and 
$\C$-anti-linear.

We can rescale the triple $(Z_\psi,H_\psi,Z_{-\psi})$   so that 
$\sigma(H_\psi)=-H_\psi$, $\sigma(Z_\psi)=-Z_{-\psi}$, 
$\sigma(Z_{-\psi})=-Z_\psi$, and the bracket 
relations  
$[Z_\psi,Z_{-\psi}]=H_\psi$, $[H_\psi,Z_\psi]=2Z_\psi$ and 
$[H_\psi,Z_{-\psi}]=-2Z_{-\psi}$   are
preserved. Let $(\cdot,\cdot)_\g$  be the  complex Killing form on $\g$ rescaled  
so that $(Z_{\psi},Z_{-\psi})_\g=\half$.
We have a positive definite Hermitian inner product
$(\cdot|\cdot)$  on $\g$ defined by $(u|v)=-(u,\sigma(v))_\g$.

Let $\g_k$ be the $k$-eigenspace of $\ad\,h$ on $\g$ where we set
\begin{equation}\label{eq:h}
h=H_\psi
\end{equation}
Since  $h$ is semisimple, $\g$ decomposes into the direct sum of the
eigenspaces $\g_k$. It is well known that the eigenvalues of $\ad\,H_{\psi}$ lie
in $\{\pm 2,\pm 1,0\}$ and $\g_{\pm 2}=\C Z_{\pm\psi}$. Thus we get the
decomposition
\begin{equation}\label{eq:g-sum}
\g=\g_2\oplus\g_1\oplus\g_0\oplus\g_{-1}\oplus\g_{-2} 
\end{equation}
We have $\dimC\g_k=\dimC\g_{-k}$ and 
the spaces $\g_{\pm 1}$ are even dimensional. We put  
$m=\half\dimC\g_{\pm 1}$. Then (see, e.g., \cite{bkLagr})
\begin{equation}\label{eq:dimO}
\dimC\O=2m+2
\end{equation}

We  now put
\begin{equation}
\label{eq:ehoe}
x_\psi=Z_\psi,\qquad x_0'=H_{\psi},\qquad x_0=Z_{-\psi}
\end{equation}
Then 
\begin{equation}
(x_0|x_0)=1/2 \AND [x'_0,x_0]=-2x_0
\end{equation}
By  \cite{abr2},
we can find a complex basis $x_1,\dots,x_m,x'_1,\dots x'_m$ of $\g_{-1}$ 
such  that
\begin{equation}
\label{eq:xbasis}
\begin{array}{c}
(x_i|x_i)=(x'_i|x'_i)=1/2\\[10pt]
[x_i,x_j]=[x'_i,x'_j]=0\AND [x'_i,x_j]=\delta_{ij}x_0 
\end{array}
\end{equation}

We have, as in (\ref{eq:fv}),   a $G $-linear  isomorphism
\begin{equation}
\g\to R_1(\O),\qquad x\mapsto f_x
\end{equation}
where $f_x(y)=(x,y)_\g$.  Also as in (\ref{eq:2-eta}) 
there is an infinitesimal vector field action 
\begin{equation}\g\to{\Vect}\,\O, \qquad  x\to\eta^x\end{equation}  
defined by differentiating the $G$-action on $\O$. 
Then $\eta^x(f_v)=f_{[x,v]}$.
We have the complex  algebra gradings 
\begin{equation}\label{eq:O-grad}
R(\O)=\bigoplus_{p\in\Z_+}R_p(\O)\AND \D(\O)=\bigoplus_{p\in\Z}\D_p(\O)
\end{equation}
defined by the Euler $\C^*$-action as in (\ref{eq:2-grad}), since
$R_p(\O)=0$ for $p<0$.

All our constructions  on $\O$ thus far have  been algebraic holomorphic. However now 
we introduce complex conjugation on $\O$ into the picture.
(See e.g. \cite[Appendix, A.5]{brInstantons}  for a discussion of complex
conjugation on varieties.)

The map $\sigma:\g\to\g$ is  a $\C$-anti-linear involution of $\g$ preserving 
$\O$.  It follows that $\sigma$, as well as $-\sigma$, is an anti-holomorphic involution  
of $\O$ which  is  also real algebraic. 
We define {\it complex conjugation} on $\O$ to be the $U$-invariant map
\begin{equation}\label{eq:cc}
\nu:\O\to\O,\qquad  \nu(z)=-\sigma(z)
\end{equation}

Now $\nu$ defines complex conjugation maps 
\begin{equation}
R(\O)\to R(\O),\quad f\mapsto\of\AND
\D(\O)\to \D(\O),\quad D\mapsto\oD
\end{equation}
 where $\of$ and $\oD$ are given  by
$\of(z)=\ovl{f(\nu(z))}$ and $\oD(f)=\ovl{D(\of)}$.
These complex conjugation maps are $\C$-anti-linear $\R$-algebra 
homomorphisms; in   particular
$\ovl{f_1f_2}=\of_1\of_2$ and $\ovl{D_1D_2}=\oD_1\oD_2$.

\begin{lem}
\label{lem:cc}
\item{\rm(i)} We have $\ovl{E}=E$  and consequently complex conjugation 
preserves the Euler gradings of $R(\O)$ and $\D(\O)$ in (\ref{eq:O-grad}).
\item{\rm(ii)} For $x\in\g$, we have 
$\ovl{f_x}=-f_{\sigma(x)}$ and $\ovl{\eta^x}=\eta^{\sigma(x)}$.
\end{lem}
{\bf Proof:} (i) is clear since $\nu$ is $\C$-anti-linear.
For (ii) we note that  
$[x,y]=[\sigma(x),\sigma(y)]$ and 
$(x,y)_\g=\ovl{(\sigma(x),\sigma(y))_\g}$ for all $x,y\in\g$
since $\sigma$ is a  anti-complex Lie algebra involution of $\g$. Consequently
\[\ovl{f_x}(z)=\ovl{f_x(-\sigma(z))}=-\ovl{(x,\sigma(z))_\g}
=-(\sigma(x),z)_\g=-f_{\sigma(x)}(z)\]
This computes  $\ovl{f_x}$. Using this we find 
\[\ovl{\eta^x}(f_y)=-\ovl{\eta^x(f_{\sigma(y)})}=
-\ovl{f_{[x,\sigma(y)]}}=f_{[\sigma(x),y]}=
\eta^{\sigma(x)}(f_y)   \]
But any vector field on $\O$ is uniquely determined by its effect on the functions     
$f_y\in R_1(\O)$. So  this proves that $\ovl{\eta^x}=\eta^{\sigma(x)}$.
\Qed

\section{Symbols and Quantization on $\O$}
\label{sec:quant}

\subsection{Set-up}
\label{ssec:setup} 

If $X$ is a complex algebraic manifold, then the  holomorphic cotangent bundle  
$T^*X\to X$  is again a complex algebraic manifold. 
The order $d$  symbol gives a linear map   $\D^d(X)\to R^d(T^*X)$,
$D\mapsto\Sigma_d(D)$, where  $R^d(T^*X)\subset R(T^*X)$ is the subspace of
functions which are homogeneous  (polynomial) functions of degree $d$   on the fibers 
of  $T^*X\to X$.  We have an algebra grading   
$R(T^*X)=\oplus_{d\in\Z_+}R^d(T^*X)$,
which we call the \emph{symbol grading}.  

If $D\in\D(X)$ has order  $d$, then $\Sigma_d(D)$ is the
\emph{principal symbol} of $D$, denoted by  $\sy\,D$.
See, e.g., \cite[Appendix]{bkDiff} for  
a resume of some basic facts about differential operators and their symbols.

An algebraic holomorphic action of $G$ on $X$ induces a natural action of $G$ on 
$T^*X$ which is  algebraic holomorphic symplectic.  Then the canonical projection
$T^*X\to X$ is $G$-equivariant; the induced action of $G$ on $T^*X$ is called the 
{\it  canonical lift} of the $G$-action on $X$.
The induced representation of  $G$  on $R(T^*X)$ is locally finite. The symbol grading
is $G$-invariant and the principal symbol map  is $G$-equivariant. 

Let  $\la\in R(T^*\O)$ be the
principal symbol of the  Euler vector field $E$ on $\O$, i.e, 
\begin{equation}\label{eq:la}
\la=\sy\, E
\end{equation}
Then $\la$ is $G$-invariant.  The canonical lift of the Euler $\C^*$-action  on $\O$ 
defines a $G$-invariant  algebra grading $R(T^*\O)=\oplus_{p\in\Z}R_p(T^*\O)$. Then
\begin{equation}\label{eq:R_p(T^*O)} 
R_p(T^*\O)=\{\phi\in R(T^*\O)\,|\, \{\la,\phi\}=p\phi\}
\end{equation}
We say that $\phi$ is \emph{Euler homogeneous of degree $p$} if 
$\phi\in  R_p(T^*\O)$.  

Now we get a  $G$-invariant complex algebra bigrading 
\begin{equation}\label{eq:bigr}
R(T^*\O)=\bigoplus_{p\in\Z, d\in\Z_+}R_{p}^d(T^*\O)
\end{equation}
where $R_{p}^d(T^*\O)=R^d(T^*\O)\cap R_{p}(T^*\O)$. The order $d$ symbol
mapping sends   $\D^d_p(\O)$ into $R_{p}^d(T^*\O)$.

Let  $\Phi^{x}\in R(T^*\O)$ be the principal symbol of the vector field $\eta^x$,  i.e., 
\begin{equation}
\label{Phi}
\Phi^{x}=\sy\, \eta^{x}\qquad\mbox{for all $x\in\g$}
\end{equation}
Then $\la,\Phi^x\in R^1_0(T^*\O)$.

For $i=0,\dots,m$ we set
\begin{equation} \label{eq:f's}
f_i=f_{x_i}\AND f'_i=f_{x'_i}\end{equation}
We recall from \cite{abr2} that these $2m+2$ functions form a set of local \'etale
coordinates on $\O$. We also put
\begin{equation}\label{eq:fpsi}
f_\psi=f_{x_\psi}
\end{equation}

We have a Zariski open dense algebraic submanifold $\Oreg$ of $\O$ defined by
\begin{equation}\label{eq:Oreg}
\Oreg=\{z\in\O\,|\, f_0(z)\neq 0\}
\end{equation}
This construction of $\Oreg$  breaks the
$G$-symmetry but not the infinitesimal $\g$-symmetry by   
vector fields. By Proposition \ref{prop:diffextend}, $R(\O)$ and $\D(\O)$
sit inside  $R(\Oreg)$  and   $\D(\Oreg)$, respectively, as the spaces of 
$\g$-finite vectors. 

We introduce the following  vector fields on $\Oreg$, where $i=1,\dots,m$,
\begin{equation}
\label{eq:Deltadef}
\Delta^{x_{i}}=\eta^{x_i}-\frac{f_i}{f_0}\eta^{x_0}\AND
\Delta^{x'_{i}}=\eta^{x'_i}-\frac{f'_i}{f_0}\eta^{x_0}
\end{equation}
and their principal symbols 
\begin{equation}
\label{eq:Theta}
\Theta^{x_{i}}=\Phi^{x_i}-\frac{f_i}{f_0}\Phi^{x_0}\AND
\Theta^{x'_{i}}=\Phi^{x'_i}-\frac{f'_i}{f_0}\Phi^{x_0}  
\end{equation}

We define a polynomial $P(w_1,w'_1,\dots,w_m,w'_m)$ by
\begin{equation}
\label{eq:Pdef}
\frac{1}{24}(\ad\,w)^4(x_{\psi})=P(w_1,w'_1,\dots,w_m,w'_m)x_0
\end{equation} 
where
\begin{equation}w=\sum_{i=1}^m(w_ix_i+w'_ix'_i)\end{equation}

\subsection{Main Results}

We assume from now on that $\g\neq\sp(2n,\C)$, $n\ge 1$. We recall
\begin{thm} \cite{abr2}
\label{thm:ES}
\item{\rm (i)}
There exists a  unique (up to scaling) non-zero $G$-linear map
\begin{equation}
\label{eq:r}
r:\g\to\RR(T^*\O),\qquad x\mapsto r_x
\end{equation}
So $r(\g)$ is the unique subspace in $\RR(T^*\O)$ which is $G$-irreducible and 
carries the adjoint representation. A  lowest weight vector $r_0=r_{x_0}$ in
$r(\g)$ is given by the formula
\begin{equation}
\label{eq:rzero}
r_0=\frac{1}{f_0}\left(P(\Theta^{x_1},\Theta^{x'_1}\dots,\Theta^{x_m},
\Theta^{x'_m})-\frac{1}{4}\la^2(\Phi^{x_{0}})^2\right)
\end{equation}  
\item{\rm (ii)} Suppose $p<4$. Then  there is no non-zero $G$-linear map     
$\g\to R_{-1}^p(T^*\O)$. I.e., $R_{-1}^p(T^*\O)$ contains no copy of the  
adjoint representation.
\end{thm}

\begin{rem} \rm
The formula (\ref{eq:rzero}) applies equally well when $\g=\sp(2n,\C)$,  $n\ge
1$. But then $P=0$. The symbol $f_0^{-1}(\Phi^{x_{0}})^2$ easily quantizes  to a
differential operator on $\O$.  See \cite{abr2}.
\end{rem}

Our main result is the $G$-equivariant quantization of these symbols $r_x$
into differential operators $D_x$ on $\O$ in  the cases where $\g$ is classical. 

\begin{thm} 
\label{thm:quant} 
Assume $\g$ is a complex simple Lie algebra  of classical type and 
$\g\neq\sp(2n,\C)$, $n\ge 1$.
Then there exists a unique differential operator $D_0$ on $\O$
such that
\item{\rm (i)} $D_0$ has order $4$ and the principal symbol of $D_0$ is $r_0$.
\item{\rm (ii)} $D_0$ is a lowest weight vector in a copy of $\g$ inside $\D(\O)$.
\item{\rm (iii)} $D_0$ is an Euler  homogeneous operator on $R(\O)$ of degree 
$-1$, i.e., $D_0\in\D_{-1}(\O)$.
\par
\noindent Moreover we construct $D_0$ as an explicit non-commutative polynomial 
in the vector fields $\Delta^{x_i},\Delta^{x'_i}$, $\eta^{x_0}$, and $E$. We prove
that
\begin{equation}\label{eq:pos}
D_0(f_\psi^k)=\gamma(k)f_\psi^{k-1}
\end{equation}
where $\gamma(k)$ is an explicit degree $4$ polynomial in $k$ with leading term
$k^4$ and all coefficients non-negative.
\end{thm} 
{\bf Proof:} The uniqueness follows from Theorem \ref{thm:ES} as it says that
$r_0$ is, up to scaling,
 the unique symbol $r\in R^{\le 4}_{-1}(T^*\O)$ such that $r$ is a lowest
weight vector of a copy of the adjoint representation. We prove the existence in the
next section by  case-by-case analysis for $\g=\sl(N,\C)$ and  $\g=\so(N,\C)$. 
See Propositions \ref{prop:SL} and \ref{prop:SO}.   \Qed
  
Property (ii) in Theorem \ref{thm:quant} implies that the $G$-subrepresentation of
$\D(\O)$ generated by
$D_0$ is irreducible and carries the   adjoint representation.
In fact  we are getting an equivariant quantization in the following sense.

\begin{cor} 
\label{cor:equiv}  There exists a  unique $G$-linear map 
\begin{equation}\label{eq: equiv}
\g\to{\D}^4_{-1}(\O),\qquad x\mapsto D_x
\end{equation}
quantizing (\ref{eq:r}) in the sense that, for all $x\in\g$,  
\begin{equation}\label{eq:sy-D}
\sy\,D_x=r_x
\end{equation}
I.e., $D_x$ has order $4$ and the principal symbol of $D_x$ is $r_x$.
Then $D_{x_0}$ is the operator  $D_0$ constructed in Theorem \ref{thm:quant}.
\end{cor}

We get two more corollaries to  Theorem \ref{thm:quant}   because of
 Propositions \ref{prop:Vcomm} and  \ref{prop:max}. The first is immediate.

\begin{cor} 
\label{cor:commute}
The differential operators $D_x$, $x\in\g$,  all commute and generate a maximal
commutative   subalgebra $\A$ of $\D(\O)$. 
Consequently  $\A=\oplus_{p\in\Z_+}\A_{-p}$ is graded in non-positive degrees 
with  $\A_0=\C$ and $\A_{-1}=\{D_x\,|\,x\in\g\}$.

$\A$ is  isomorphic to $R(\O)$. In fact, the  mapping $f_x\mapsto D_x$ extends 
uniquely to a  $G$-equivariant complex algebra isomorphism
\begin{equation}  R(\O)\to\A,\qquad f\mapsto D_f  \end{equation} 
\end{cor}

\begin{cor} 
\label{cor:Herm}
$R(\O)$ admits a unique positive-definite  Hermitian inner product $(\cdot|\cdot)$ 
invariant under the maximal compact subgroup $U$ of $G$ such
that $(1|1)=1$ and  the operators $\of$ and $D_f$ are adjoint for all  $f\in
R(\O)$. Then the formula (\ref{eq:HIP}) holds. 
\end{cor} 
{\bf Proof:}
Proposition \ref{prop:max} applies because   (\ref{eq:pos}) says that $D_0$ is
non-zero on $R_p(X)$ for $p\ge 1$. This  gives us a 
unique non-degenerate
$U$-invariant Hermitian inner product $(\cdot|\cdot)$ on $R(\O)$ such that
$(1|1)=1$ and    $\of$ and $D_f$ are adjoint. 

We need to check  that $(\cdot|\cdot)$  is  positive-definite  on each space
$R_p(\O)$, $p\ge 1$.   It suffices to check that $(f_\psi^p|f_\psi^p)$ is 
positive  for $p\ge 1$.  Lemma \ref{lem:cc}(ii) implies that 
\begin{equation}\label{eq:fpsibar}
\ovl{f_\psi}=f_0 
\end{equation}
since $-\sigma(x_{\psi})=x_0$ (see Subsection \ref{ssec:setup}).   
So (\ref{eq:HIP}) and (\ref{eq:pos})  give 
\begin{equation}
(f_\psi^p|f_\psi^p)=D_{f_0^p}(f_\psi^p)=
D_0^p(f_\psi^p)=\gamma(1)\cdots\gamma(p)
\end{equation}
which is positive by Theorem \ref{thm:quant}.  \Qed

\begin{rem}\rm
We will show in \cite{abr3} that $\ovl{D_x}=-D_{\sigma(x)}$. 
\end{rem}

\subsection{Strategy for quantizing  the  symbol $r_0$}
\label{ssec:strat}

Here is our strategy  for quantizing the symbol $r_0$ into the operator $D_0$.  We  
start from  the formula   (\ref{eq:rzero}) for the   symbol $r_0$.  This says that 
$S=f_0r_0$  where
\begin{equation}
\label{eq:S}
S=P(\Theta^{x_1},\Theta^{x'_1}\dots,\Theta^{x_m},\Theta^{x'_m})-
\frac{1}{4}\la^2(\Phi^{x_{0}})^2
\end{equation} 
Our idea is to construct $D_0$ by first constructing a    suitable quantization 
$\bS$ of $S$ which is left divisible by $f_0$,  and  then putting 
$D_0=f_0\i\bS$.

Now $S$ is an explicit  polynomial (once we compute $P$)  in the symbols
$\Th^{x_i},\Th^{x'_i},\la,\Phi^{x_0}$ of the vector fields 
$\Delta^{x_i},\Delta^{x'_i},E,\eta^{x_0}$ on $\Oreg$.  Hence $S$ belongs to
$R(T^*\Oreg)$.  In fact by Theorem
\ref{thm:ES}, we know that $S$ belongs to $R(T^*\O)$.

The first step of our strategy is to quantize $S$ into a differential operator $\bS$  
on $\Oreg$ using the  philosophy of  symmetrization from Weyl quantization. 
This means that we obtain  $\bS$ by taking the
expression for $S$, replacing each symbol $\Th^{x_i},\Th^{x'_i},\la,\Phi^{x_0}$ 
by the corresponding vector field, and symmetrizing to allow for any order
ambiguities. E.g.,
$\Th^{x'_i}\Th^{x_i}$ quantizes to 
$\half(\Delta^{x'_i}\Delta^{x_i}+\Delta^{x_i}\Delta^{x'_i})$.
We will refer to this procedure simply as ``Weyl quantization".
(This  does not determine the quantization uniquely, as one can symmetrize the 
symbol $(\Th^{x'_i}\Th^{x_i})^2$ in different ways;  we merely choose a 
convenient scheme.) 

Thus we construct $\bS$ as a   non-commutative polynomial    in
the vector fields  $\Delta^{x_i},\Delta^{x'_i},E,\eta^{x_0}$ on $\Oreg$.  
We want to quantize $S$ into $\bS$ in a $G$-equivariant way
so that $\bS$, like $S$, is a lowest weight vector of a copy of 
$V_{2\psi}\simeq\g^{\boxtimes 2}$. It turns out to
be too hard (and unnecessary) to verify all this right away, but we do construct  
$\bS$ so that the following properties are evident:
\begin{equation}
\label{eq:Sreq}
[\eta^z,\bS]=0\mbox{ for all $z\in\heis$ \quad and }\quad  [\eta^h,\bS]=-4\bS
\end{equation}
where 
\begin{equation}
\label{eq:heis}
\heis=\g_{-1}\oplus\g_{-2}
\end{equation}

The conditions in (\ref{eq:Sreq}) go half-way towards proving that
$\bS$  is a lowest weight vector of a copy of $\g^{\boxtimes 2}$.
Indeed,  we have the  following easy fact from the  Cartan-Weyl theory of 
representations.

\begin{lem}
\label{lem:rep}
A vector $v$ in a  representation $\rho:\g\to\End\, W$
is a lowest weight vector of  a copy of $V_{k\psi}\simeq\g^{\boxtimes k}$
if and only if  the following four conditions are satisfied:
\item{\rm (i)}  $v$ is $\g$-finite
\item{\rm (ii)} $\rho_h(v)=-2kv$
\item{\rm (iii)} $\rho_x(v)=0$ for all $x\in\g_0'$ where  $\g_0'$ is the  
orthogonal  space to $h$ in $\g_0$ 
\item{\rm (iv)} $\rho_z(v)=0$ for all $z\in\heis$
\end{lem}

The second  step is to show that $f_0\i\bS$ belongs to $\D(\O)$,
(We  skip the intermediate step of showing that $\bS$ belongs to
$\D(\O)$.)   In fact, the multiplication operator $f_0$  commutes with all the vector
fields $\Delta^{x_i},\Delta^{x'_i},E,\eta^{x_0}$, and so commutes with $\bS$.
Thus we  need not distinguish between $f_0\i\bS$ and $\bS f_0\i$.

We accomplish the second step by implementing Proposition \ref{prop:extend}
below. Before giving  that result, we first prove two  more general Lemmas
about  extending differential operators from $\Oreg$ to $\O$.

\begin{lem}
\label{lem:extendgen}
Suppose $\bT$ is a differential operator on $\Oreg$. Then  $\bT$ extends to a 
differential operator on
$\O$ if and only if  the operator $\bT:R(\Oreg)\to R(\Oreg)$ satisfies
$\bT(R(\O))\subset R(\O)$.
\end{lem}
{\bf Proof:} Immediate as $\O$ is quasi-affine; see e.g., 
\cite[Appendix, A.4 and A.6]{bkDiff}.
\Qed

\begin{lem}
\label{lem:extendlow}
Suppose $\bT$ is a differential operator on $\Oreg$ such that $[\eta^z,\bT]=0$ 
for all
$z\in\heis$.  Then  $\bT$ extends to a differential operator on $\O$
if and only if  $\bT(f_{\psi}^k)\in R(\O)$ for $k=0,1,2\dots$. 
\end{lem}
{\bf Proof:} The functions  $f_{\psi}^k$, $k=0,1,2\dots$, form a complete set of 
highest weight vectors in $R(\O)$ for the $\g$-representation given by the vector
fields $\eta^x$. Precisely, 
$f_{\psi}^k$    generates  $R_k(\O)$ under the action of $\n^-$ so that
$\U(\n^-)\cdot f_{\psi}^k=R_k(\O)$.  We have $\n^-\subset\g_0\oplus\heis$
and also $\g_0\cdot f_{\psi}^k=\C f_{\psi}^k$ by Lemma  \ref{lem:rep}(iii).
It follows that  
\begin{equation}
\label{eq:Uheis}
R_k(\O)=\U(\heis)\cdot f_{\psi}^k
\end{equation}
Now suppose $\bT(f_{\psi}^k)\in R(\O)$. Then, since $\bT$ commutes with
$\eta^z$ for all $z\in\heis$, we   get
\[\bT(R_k(\O))=\bT(\U(\heis)\cdot f_{\psi}^k)=\U(\heis)\cdot
\bT(f_{\psi}^k)\subset R(\O)\]
Thus $\bT(R(\O))\subset R(\O)$ and  so $\bT$ extends to a differential operator 
by  Lemma \ref{lem:extendgen}. The  converse is obvious.  \Qed
  
To show  that  $f_0\i\bS$ belongs to $\D(\O)$ we will use

\begin{prop}
\label{prop:extend}
Suppose  $\bS$ is a  differential operator on $\O\reg$ such that
$\bS$ is Euler homogeneous of degree $0$ and $\bS$ satisfies (\ref{eq:Sreq}).
Then   the following two  properties are equivalent:
\item{\rm (i)} $f_0^{-1}\bS\in\D(\O)$ 
\item{\rm (ii)} $\bS(f_{\psi}^k)=\gamma_kf_0f_{\psi}^{k-1}$ for    
$k=0,1,2\dots$ where  $\gamma_0,\gamma_1,\dots$ are scalars.
\mypar
Moreover (i) and (ii) imply 
\item{\rm (iii)}
$f_0^{-1}\bS$  is a  lowest weight vector of a copy of $\g$  in $\D_{-1}(\O)$ .
\end{prop}
{\bf Proof:}
Put   $\bT=f_0\i\bS$.  Then $[\eta^z,\bT]=0$ 
for all  $z\in\heis$ since $\eta^z$ commutes with $f_0$ and $\bS$. So
Lemma \ref{lem:extendlow} applies and says that $\bT\in\D(\O)$ if and only if
$\bT(f_{\psi}^k)\in R(\O)$ for all $k=0,1,\dots$. Hence (ii) implies (i).
 
Conversely suppose (i) holds. Then $g=\bT(f_{\psi}^k)$ lies in $R_{k-1}(\O)$
and $\eta^h(g)=2(k-1)g$. But it  follows using  (\ref{eq:Uheis}) 
that $f_{\psi}^{k-1}$ is, up to scaling, the unique eigenfunction in $R_{k-1}(\O)$
of $\eta^h$  with eigenvalue $2(k-1)$. So $g=\gamma_k f_{\psi}^{k-1}$.
This proves (ii).

Finally we show that properties (i) and (ii) imply (iii).  We can apply Lemma
\ref{lem:rep} where $W=\D(\O)$, $v=\bT=f_0\i\bS$ and $k=1$. Now condition  (i) 
in Lemma \ref{lem:rep} is satisfied since the whole $\g$-representation $\D(\O)$
is locally finite. Also  
conditions (ii) and  (iv) follow from (\ref{eq:Sreq}). So we need to check (iii).

Let $x\in\g_0'$. We want to show that the differential
operator $[\eta^x,\bT]$ on $\O$ is zero.  We can do this by showing   
that  $[\eta^x,\bT]$ annihilates $R(\O)$.  Now $\eta^x(f_{\psi})=0$
(use (\ref{eq:g-sum})) and so $\bT(f_{\psi}^k)=\gamma_k f_{\psi}^{k-1}$ implies
$[\eta^x,\bT](f_{\psi}^k)=0$. 
But also $[\eta^x,\bT]$ commutes with the vector fields $\eta^z$, $z\in\heis$.
Indeed $[\eta^z, [\eta^x,\bT]]=[\eta^{[z,x]},\bT]=[\eta^{z'},\bT]=0$
as $z'=[z,x]\in\heis$. So then $[\eta^x,\bT]$ annihilates $R(\O)$
because of  (\ref{eq:Uheis}).
\Qed

In carrying out this plan for quantizing $r_0$,  we encountered a problem as we found  
no ``Weyl  quantization" of the term $\la^2(\Phi^{x_{0}})^2$ in (\ref{eq:S}) which 
produced an operator $\bS$ satisfying the condition
$\bS(f_{\psi}^k)=\gamma_kf_0f_{\psi}^{k-1}$.    (Such a Weyl quantization may or 
may  not be possible here; our ``Weyl quantization" of the first term is already not
unique.  Probably Weyl quantization fails to give the operators we want; cf. the 
quantization in \cite{lo}.)

To get around this, we applied   Weyl quantization to the first term in  
(\ref{eq:S}) after computing $P$ explicitly.  Then we ``guessed" that 
$\la^2(\Phi^{x_{0}})^2$ should quantize into the operator
\begin{equation}
\label{eq:EE}
(E+c_1)(E+c_2)(\eta^{x_0})^2
\end{equation}
where $c_1$ and $c_2$ were unknown constants.  There is no problem in this as
\begin{lem}
\label{lem:qeta}
For any  polynomial $q(E)$ in the Euler vector field, the operator
$q(E)(\eta^{x_0})^2$ is a lowest weight vector of a copy of  
$V_{2\psi}\simeq\g^{\boxtimes 2}$ in $\D(\O)$.
\end{lem}

Thus we got  a candidate for $\bS$ with
unknown parameters $c_1$ and $c_2$. We then calculated $\bS(f_{\psi}^k)$ and 
``solved" for $c_1$ and $c_2$ by demanding that $\bS(f_{\psi}^k)$ be a multiple 
of  $f_0f_{\psi}^{k-1}$. There was no guarantee  a priori that solutions for  
$c_1$ and $c_2$ existed, but in fact we  found a unique solution for the polynomial
$(E+c_1)(E+c_2)$.

In our presentation below, we  in fact use a general polynomial $q(E)$ in place
of the particular polynomial $(E+c_1)(E+c_2)$. We show that the required relation
$\bS(f_{\psi}^k)=\gamma_kf_0f_{\psi}^{k-1}$ forces $q(E)$ to be of the form
$(E+c_1)(E+c_2)$ and we determine $c_1$ and $c_2$.

In working this out,  we do not yet have a unified treatment of all 
complex simple Lie algebras. Instead we  treat only the cases where  $\g$ is
classical.  We  work  out the polynomial $P$ and  the operator $\bS$ in Section
\ref{sec:construct}, first for   $\g=\sl(N,\C)$ and then for  $\g=\so(N,\C)$. 

We conjecture that Theorem \ref{thm:quant} holds also for the five exceptional
simple  complex Lie algebras.

\section{Explicit construction of  $D_0$ for $\g$  classical}
\label{sec:construct}

\subsection{Coordinate representations of vector fields}
\label{ssec:4.1} 
 
Our aim is to construct $D_0$ as the quantization of  the symbol $r_0$.
Our   strategy,   developed in  Subsection \ref{ssec:strat},  is to quantize the symbol
$S=f_0r_0$ into an operator $\bS$ which satisfies condition (ii) in 
Prop \ref{prop:extend}. To construct $\bS$ and compute $\bS(f_{\psi}^k)$,
we will simply work out everything in terms of  our local coordinates 
$f_0,f'_0,f_i,f'_i$  on $\O$ from (\ref{eq:f's}).
 
Fortunately, in  proving Theorem \ref{thm:ES} in \cite{abr2},
we first found a formula for the highest weight vector
$f_{\psi}\in R_1(\O)\simeq\g$. This was:
\begin{prop}{\bf \cite{abr2}} 
\label{prop:fe}
The unique expression for the function $f_{\psi}\in R_1(\O)$
in  terms of our local coordinates  $f_0,f'_0,f_i,f'_i$, $i=1,\dots,m$, on $\O$ is 
\begin{equation}
\label{eq:fe}
f_\psi=\frac{1}{f_0^3}P(f_1,f_1',\dots,f_1,f_1')-\frac{1}{4}\frac{(f_0')^2}{f_0}
\end{equation}
\end{prop}

The unique expressions for  the vector fields $\eta^{x_i},\eta^{x'_i},\eta^{x_0}$ 
in terms of our local coordinates  $f_i,f'_i,f_0,f'_0$ are, where $i=1,\dots,m$, 
\begin{equation}
\label{eq:etacalc}
\eta^{x_i}=-f_0\pdb{f'_i}+f_i\pdb{f'_{0}}\qquad
\eta^{x'_i}=f_0\pdb{f_i}+f'_i\pdb{f'_{0}},\qquad
\eta^{x_{0}}=2f_0\pdb{f'_0}  
\end{equation}
These follow as    $\eta^x(f_y)=f_{[x,y]}$ for any $x,y\in\g$. Then
(\ref{eq:Deltadef}) gives
\begin{equation}
\label{eq:Deltacalc}
\Delta^{x_{i}}=-f_0\pdb{f'_i}-f_i\pdb{f'_{0}} ,\qquad 
\Delta^{x'_{i}}=f_0\pdb{f_i}-f'_i\pdb{f'_{0}}
\end{equation}

The operators $\Delta^{x_i},\Delta^{x'_i},\eta^{x_0}$ also span a Heisenberg 
Lie algebra  since 
\begin{equation}
\label{eq:Deltabr}
[\Delta^{x'_{i}},\Delta^{x'_{j}}]=[\Delta^{x_{i}},\Delta^{x_{j}}]
=[\Delta^{x'_{i}},\eta^{x_0}]=[\Delta^{x_{i}},\eta^{x_0}]=0,\qquad 
[\Delta^{x_{i}},\Delta^{x'_{j}}]=2\delta_{ij}\eta^{x_0}
\end{equation}

The next two facts are trivial to verify.
\begin{lem}
\label{lem:Deltacomm}
The operators $\Delta^{x_i},\Delta^{x'_i}$ all commute with the vector fields 
$\eta^z$ where $z$ lies in the Heisenberg Lie algebra $\heis$.
\end{lem}

\begin{lem}
\label{lem:etah}
We have the relations $\eta^h(f_0)=-2f_0$, $\eta^h(f'_0)=0$, 
$\eta^h(f_i)=-f_i$, $\eta^h(f'_i)=-f'_i$  for $i=1,\dots,m$, and 
$[\eta^h,\Delta^y]=-\Delta^y$ if $y\in\g_{-1}$.
\end{lem}

Also, we have $\eta^{x_0}f_\psi=f_{[x_0,x_\psi]}=f_{-x_0'}=-f_0'$ and so using
(\ref{eq:etacalc})  we find
\begin{equation}
\label{eq:etasq}
(\eta^{x_0})^2f_\psi^k=
\left[-2kf_0f_\psi+ k(k-1)(f_0')^2\right]f_\psi^{k-2}
\end{equation}

To reduce the number of formulas, we introduce the superscript $\ep$ 
with $\ep=\pm 1$ where $\ep=1$ indicates primed quantities while $\ep=-1$ 
indicates unprimed quantities. Thus $x_i^\ep=x'_i$  if $\ep=1$ while 
$x_i^\ep=x_i$ if   $\ep=-1$ and so on. Then
(\ref{eq:etacalc}) and (\ref{eq:Deltacalc}) give
 \begin{equation}
\label{eq:etaDelta}
\eta^{x_i^\ep}=\ep f_0\pdb{f^{-\ep}_i}+f^\ep_i\pdb{f'_{0}}\qquad
\Delta^{x_{i}^\ep}=\ep f_0\pdb{f^{-\ep}_i}-f^\ep_i\pdb{f'_{0}}   
\end{equation}

\subsection{The case $\g=\sl(n+1,\C)$, $n\ge 2$}

In this case $\g$ is of type $A_n$. The minimal nilpotent orbit $\O$ is
\begin{equation}
\label{eq:OSL}
\O=\{X\in\g\,|\, X^2=0, \mbox{rank}(X)=1\}
\end{equation}
Then  $\dimC\O=2n$  and so 
\begin{equation} 
\label{eq:mSL}
m=n-1  
\end{equation}
A   standard basis of $\gl(n+1,\C)$ is given by the elementary matrices 
$E_{i,j}$, 
$0\le i,j\le n$. We   choose $x_\psi=E_{0,n}$, $x_0'=E_{0,0}-E_{n,n}$, and
$x_0=E_{n,0}$. 
Then the normalized Killing form  on $\g$ is given by  
$(X,X')_\g=\half(\mbox{Trace}\, XX')$.  We take the  Cartan involution 
$\sigma$ to  be  $\sigma(X)=-X^*$.

The matrices $x'_p=E_{n,p}$ and $x_p=E_{p,0}$,
$p=1,\dots,m$,  form a basis of $\g_{-1}$ and  satisfy the conditions  in
(\ref{eq:xbasis}).  So we get
\begin{equation}
\label{eq:wmatrixSL}
w=\sum_{p=1}^m(w_pE_{p,0}+w_p'E_{n,p})=\left(\begin{array}{ccccc}   
0&0&\dots&0&0\\ w_1&0&\dots&0&0\\
\vdots&\vdots&&\vdots&\vdots\\ w_m&0&\dots&0&0\\
0& w'_1&\dots&  w'_m& 0\end{array}\right)
\end{equation}
 
Next we need to compute the polynomial $P$ defined in (\ref{eq:Pdef}).
We have the standard  operator identity
$\mbox{Ad}(\exp tw)=\sum_{k=0}^\infty\frac{1}{k!}(\ad\,tw)^k$
and so (\ref{eq:Pdef}) gives
\begin{equation}
\label{eq:Pcoeff}
\mbox{$P(w_i,w'_i)$ is the coefficient of $t^4x_0$ in  
$(\mbox{Ad}(\exp tw))\cdot x_\psi$}
\end{equation}
We can easily compute this since  
\begin{equation}
\label{eq:Ad}
(\mbox{Ad}(\exp tw))\cdot x_\psi=(\exp tw)x_\psi(\exp -tw)  
\end{equation}
where the right side is a matrix product. We find
\begin{equation}
\label{eq:exptwSL}
\exp(tw)=I+tw+\frac{1}{2}t^2
\textstyle\left(\sum_{p=1}^mw_pw'_p\right)E_{n,0}
\end{equation}
since   $w^2=(\sum_{p=1}^mw_pw'_p)E_{n,0}$ and   $w^3=0$. Then we get
\begin{equation}
\label{eq:PSL}
P(w_1,w_1',\dots,w_m,w_m')=\frac{1}{4}
\textstyle\left(\sum_{p=1}^mw_pw'_p\right)^2
\end{equation}

Now plugging this expression for $P$ into  (\ref{eq:fe})  and(\ref{eq:rzero})  
we get 
\begin{equation}
\label{eq:ferSL}
f_\psi=\frac{1}{4f_0^3}\left[a^2-(f_0f'_0)^2\right]\AND
r_0=\frac{1}{4f_0}\left[A^2-\la^2(\Phi^{x_0})^2\right] 
\end{equation}
where
\begin{equation}
\label{eq:aASL}
a=\sum_{i=1}^mf_if'_i\AND A=\sum_{i=1}^m\Th^{x_i}\Th^{x'_i}
\end{equation}
 
To begin the process of quantizing $S=f_0r_0$, we  quantize $A$ 
into the operator
\begin{equation}
\label{eq:bASL}
\bA=\half \sum_{i=1}^m(\Delta^{x_i}\Delta^{x'_i}+\Delta^{x'_i}\Delta^{x_i})
\end{equation}

\begin{lem}
\label{lem:ASL}
$\bA$  is a  differential operator on $\O\reg$ of order  $2$ with  principal symbol  
$A$. $\bA$ is Euler homogeneous of degree $0$.
$\bA$ satisfies $[\eta^z,\bA]=0$ for all $z\in\heis$ and $[\eta^h,\bA]=-2\bA$.
\end{lem}
{\bf Proof:} The first statement is clear.  The commutativity $[\eta^z,\bA]=0$ 
follows by   Lemma \ref{lem:Deltacomm}  and the weight relation
$[\eta^h,\bA]=-2\bA$ follows by Lemma \ref{lem:etah}. 
\Qed

By calculation, we find the following two formulas, where the answer is expressed  
in  terms of our  ``$\ep$" notation defined at the end of   Subsection 
\ref{ssec:4.1}.  
\begin{equation}
\label{eq:DeltaaSL} 
\Delta^{x_i^\ep}(a)=\ep f_0f_i^\ep\AND   \Delta^{x_i^\ep}(f_0f'_0)=
-f_0f_i^\ep
\end{equation}
So we get
\begin{equation} 
\label{eq:DeltafeSL}
\Delta^{x_i^\ep}(f_\psi)=\frac{\ep}{2}\frac{f_i^\ep(a+\ep f_0f_0')}{f_0^2}
\end{equation}

Next we determine $\bA(f_\psi^k)$ starting from the identity
\begin{equation}
\label{eq:bAfekgenSL}
\bA(f_{\psi}^k)=k(\bA f_{\psi})f_{\psi}^{k-1}+k(k-1)f_{\psi}^{k-2}\sum_{i=1}^m
(\Delta^{x_i}f_{\psi})(\Delta^{x'_i}f_{\psi})
\end{equation}
Using (\ref{eq:DeltafeSL}) and  (\ref{eq:ferSL})  we find
 $\Delta^{x_i}(f_{\psi})\Delta^{x'_i}(f_{\psi})=-f_if_i'f_{\psi}/f_0$ and so
\begin{equation}
\sum_{i=1}^m(\Delta^{x_i}f_{\psi})(\Delta^{x'_i}f_{\psi})=-\frac{af_{\psi}}{f_0}
\end{equation}
Now we need the first term in (\ref{eq:bAfekgenSL}). First  
(\ref{eq:DeltafeSL})  and  (\ref{eq:DeltaaSL}) give
\begin{equation}
\Delta^{x_i^\ep}\Delta^{x_i^{-\ep}}(f_{\psi})=\frac{-a+\ep f_0f_0'-2f_if_i'}
{2f_0}
\end{equation}
and consequently
\begin{equation}
\half(\Delta^{x_i}\Delta^{x'_i}+\Delta^{x'_i}\Delta^{x_i})(f_{\psi})=
-\frac{a+2f_if'_i}{2f_0}
\end{equation}
Now summing over $i=1,\dots,m$ we get
\begin{equation}
\label{eq:bAfeSL}
\bA(f_{\psi})=-{\textstyle\left(1+\frac{m}{2}\right)}\frac{a}{f_0}
\end{equation}
Now substituting into (\ref{eq:bAfekgenSL}) we get
\begin{equation}
\label{eq:bAfekSL}
\bA(f_{\psi}^k)=-{\textstyle k(k+\frac{m}{2})}\frac{af_{\psi}^{k-1}}{f_0}
\end{equation}

We need to  calculate $\bA^2(f_{\psi}^k)$. Since $\bA$ commutes with 
multiplication by any power of $f_0$, (\ref{eq:bAfekSL}) gives

\begin{equation}\label{eq:newA2SL} 
\bA^2(f_{\psi}^k)=-k\left(k+{\textstyle\frac{m}{2}}\right) f_0\i
\bA(af_{\psi}^{k-1})
\end{equation}
So the problem is to compute
\begin{eqnarray}
\bA(af_{\psi}^{k-1})&=&(\bA a)f_{\psi}^{k-1}+a(\bA f_{\psi}^{k-1})
\nonumber\\[4pt]
&&+(k-1)f_{\psi}^{k-2}
\sum_{i=1}^m(\Delta^{x_i}a)(\Delta^{x'_i}f_{\psi})+(\Delta^{x_i}f_{\psi})
(\Delta^{x'_i}a)
\end{eqnarray}

Easily $\bA(a)=-mf_0^2$ and also (\ref{eq:DeltafeSL}) and  
(\ref{eq:DeltaaSL})  give
\begin{equation}
(\Delta^{x_i^{\ep}}a)(\Delta^{x_i^{-\ep}}f_{\psi})=
\frac{f_if'_i}{2f_0}(-a+\ep f_0f_0')
\end{equation}
So
\begin{equation}
\sum_{i=1}^m(\Delta^{x_i}a)(\Delta^{x'_i}f_{\psi})+
(\Delta^{x_i}f_{\psi})(\Delta^{x'_i}a)=
-\frac{a^2}{f_0}
\end{equation}
Collecting all these terms we find
\begin{equation}
\bA(af_{\psi}^{k-1})=-mf_0^2f_{\psi}^{k-1}-
(k-1)\left(k+{\textstyle\frac{m}{2}}\right)\frac{a^2}{f_0}f_{\psi}^{k-2}
\end{equation}
Thus   (\ref{eq:newA2SL}) gives 
\begin{equation}
\label{eq:bAsSL}
\bA^2(f_{\psi}^k)=\left[mk\left(k+{\textstyle\frac{m}{2}}\right)f_0f_{\psi}
+k(k-1)(k+{\textstyle\frac{m}{2}})^2\frac{a^2}{f_0^2}\right]f_{\psi}^{k-2}
\end{equation}
 
Our aim is to find a polynomial $q(E)$ in the Euler vector field $E$ such that
$\bS=\frac{1}{4}\left(\bA^2-q(E)(\eta^{x_0})^2\right)$ satisfies 
$\bS(f_{\psi}^k)=\alpha_kf_0f_{\psi}^{k-1}$ for some scalars $\alpha_k$.
A priori, there is no guarantee that such an operator $q(E)$ exists.
But fortunately  (\ref{eq:ferSL}) gives the relation
\begin{equation}
\label{eq:relSL}
\frac{a^2}{f_0^2}-(f_0')^2=4f_0f_{\psi}
\end{equation}
This is the unique linear relation in $R(\Oreg)$ between  the 
three functions $\frac{a^2}{f_0^2}$, $(f_0')^2$ and $f_0f_{\psi}$.
We have $q(E)(f_{\psi}^k)=q(k)f_{\psi}^k$.
Now comparing (\ref{eq:bAsSL}) with (\ref{eq:etasq}) we  
conclude that $q(E)$ exists if and only if 
$k(k-1)(k+\frac{m}{2})^2=q(k)k(k-1)$.
Thus $q(k)=(k+\frac{m}{2})^2$ is the unique solution. 
Then arithmetic gives a formula for $\alpha_k$ (which turns out to simplify very
nicely). This proves

\begin{lem} 
\label{lem:mainSL}
Let $\bS=\frac{1}{4}\left(\bA^2-q(E)(\eta^{x_0})^2\right)$ 
where $q(E)$ is some polynomial in the Euler operator $E$.
Then $\bS$ is a differential operator on $\O\reg$ 
of order $4$ with principal symbol  $S$.  $\bS$ is Euler homogeneous of degree 
$0$ and satisfies the conditions in  (\ref{eq:Sreq}).  The  polynomial  
\begin{equation}
\label{eq:qSL}
q(E)=\textstyle\left(E+\frac{m}{2}\right)^2
\end{equation} 
is the unique polynomial in $E$ such that  $\bS$ satisfies (ii) in 
Proposition \ref{prop:extend}. Then $\bS(f_{\psi}^{k})=\alpha_kf_0f_{\psi}^{k-1}$ 
where
\begin{equation}
\label{eq:alphaSL}
\alpha_k=k^2(k+\textstyle\frac{m-1}{2})(k+\frac{m}{2})
\end{equation}
\end{lem} 

The polynomial $\alpha_k$  has leading term $k^4$ and all coefficients non-negative
since $m\ge 1$. Proposition \ref{prop:extend} now says that
$f_0^{-1}\bS$  is a  lowest weight vector of a copy of $\g$  in $\D_{-1}(\O)$ 
Thus we have proven

\begin{prop}
\label{prop:SL} 
Let $\g=\sl(m+1,\C)$, $m\ge 1$, and define $\bA$ by (\ref{eq:bASL}).
The differential operator
\begin{equation}
\label{eq:DSL}
D_0=\frac{1}{f_0}\bS=
\frac{1}{4f_0}\left[\bA^2-\textstyle(E+\frac{m}{2})^2 
(\eta^{x_0})^2\right]
\end{equation} 
on $\Oreg$ extends to an algebraic differential operator  on $\O$.
$D_0$ satisfies  all three conditions in Theorem \ref{thm:quant}.  We have
\begin{equation}
\label{eq:DkSL}
D_0(f_{\psi}^{k})=\alpha_kf_{\psi}^{k-1}
\end{equation} 
where  $\alpha_k$ is given by  (\ref{eq:alphaSL}).
\end{prop} 
 
\begin{rem}\rm
We can show that  $\bA$ extends to a differential operator on $\O$.
\end{rem}

\subsection{The case $\g=\so(N,\C)$, $N\ge 6$}

Here $\g$ is of type $D_n$ if $N=2n$ is even  or of type  $B_n$ if 
$N=2n+1$ is odd. We unify our treatment of the two types
by using the obvious inclusion $\so(2n,\C)\subset\so(2n+1,\C)$. 
The ``final answer"  given in 
Proposition  \ref{prop:SO} below is stated in a uniform way. 

To begin with, we introduce a convenient model for  the complex Lie
algebra $\g=\so(N,\C)$.   Let $u\cdot v$  be a symmetric non-degenerate
complex bilinear form  on $\C^N$. We have a complex linear Lie bracket on
$\wedge^2\C^N$  given by
\begin{equation}[a\wedge b,c\wedge d]=
(a\cdot c)b\wedge d+(b\cdot d)a\wedge c-(a\cdot d)b\wedge c-
(b\cdot c)a\wedge d
\end{equation}
Each vector $u\wedge v$ defines a skew-symmetric linear transformation
$L_{u\wedge v}$ on  $\C^n$ by
\begin{equation}
L_{u\wedge v}(a)=(u\cdot a)v-(v\cdot a)u
\end{equation} 
Extending linearly, we obtain a  
natural complex Lie algebra isomorphism  $\wedge^2\C^N\to\so(N,\C)$, 
$z\mapsto L_z$. We use this to identify $\g$ with $\wedge^2\C^N$.

Now the minimal nilpotent orbit $\O$ is   
\begin{equation}
\label{eq:OSO}
\O=\{u\wedge v\in\g=\wdg^2\C^N\,|\, u\cdot u=u\cdot v=v\cdot v=0,\quad
u\wedge v\neq 0\}
\end{equation}
The complex dimension of  $\O$ is $\dimC\O=2N-6$, and so 
\begin{equation}
\label{eq:mSO}
m=N-4
\end{equation}
Thus $m=2n-4$ if $\g$ has type $D_n$, while  $m=2n-3$ if $\g$ has type 
$B_n$.
 
We next fix a convenient basis of $\C^N$.
If $N=2n+1$, we let $v_0$ denote a  vector in $\C^N$ such that 
$v_0\cdot v_0=1$;  if $N=2n$,  we set $v_0=0$.  
Then we have a  direct sum decomposition $\C^N=V\oplus\C v_0$ where $V$ is the
orthogonal space to $v_0$.  We fix a  basis  $v_1,v'_1,\dots,v_n,v'_n$ of $V$  
such  that
$v_i\cdot v_j=v'_i\cdot v'_j=0$, and $v_i\cdot v_j'=\delta_{ij}$.
This basis of $V$, together with $v_0$ if $N$ is odd,  is a basis of $\C^N$.

Now we choose  
\begin{equation}
\label{eq:tripleSO} 
x_0=-v_{n-1}\wedge v_n,\qquad x_\psi=v'_{n-1}\wedge v'_n,\qquad 
x_0'=v_{n-1}\wedge v'_{n-1}+ v_{n}\wedge v'_n
\end{equation}
The normalized Killing form on $\g$ is given by
\begin{equation}
\label{eq:formSO}
(a\wedge b,c\wedge d)_\g=\half\left[(a\cdot c)(b\cdot d)-
(a\cdot d)(b\cdot c)\right]
\end{equation}
We take the Cartan involution $\sigma$ of $\g$ 
to be the one    induced in the natural way from the $\C$-anti-linear involution 
$\sigma:\C^N\to\C^N$  such that 
$\sigma(v_i)=v_i'$ and  $\sigma(v_0)=v_0$.

We can now construct  a nice  basis of $\g_{-1}$. This basis is given by  the set 
$v_i^\ep\wdg v_{n-1},v_i^\ep\wdg v_{n}$ where
$i=1,\dots,n-2$, $\ep=\pm 1$, together with 
$v_0\wdg v_{n-1},v_0\wdg v_{n}$ if $N$ is odd. We label the basis vectors 
in the  following way:
\begin{equation}
\label{eq:txzeroSO}
\tx_0=v_0\wedge v_{n-1},\qquad \tx'_0=v_0\wedge v_{n}
\end{equation}
and for $i=1,\dots,n-2$,
\begin{equation}
\label{eq:xbasisSO}
x_i=v_i\wdg v_{n-1},\qquad \tx_i= v'_i\wdg v_{n-1},\qquad 
x'_i=v'_i\wdg v_{n},\qquad \tx'_i=v_i\wdg v_{n}
\end{equation}
The bracket relations  among basis vectors are all zero except  for the following ones
\begin{equation}
\label{eq:bracketSO}
[x'_i,x_j]=[\tx'_i,\tx_j]=\delta_{ij}x_0,\qquad  [\tx'_0,\tx_0]=x_0
\end{equation}
where  $i,j\in\{1,\dots,n-2\}$. 

Referring back to Subsection \ref{ssec:2.3}, we see now that all the relations in 
(\ref{eq:xbasis}) are satisfied if we take 
$x_1,\dots,x_{n-2},\tx_0,\tx_1,\dots,\tx_{n-2}$  and 
$x'_1,\dots,x'_{n-2},\tx'_0,\tx'_1,\dots,\tx'_{n-2}$ to be, respectively,
the lists  $x_1,\dots,x_m$  and   $x'_1,\dots,x'_m$. 
(Here we ignore $\tx_0$ and $\tx'_0$ 
if $N$ is even.) We have broken the sets $x_1,\dots,x_m$ and $x'_1,\dots,x'_m$
into these peculiar groups because this is necessary to write down $P$ in  the next
step.

Now we put
\begin{equation}
\label{eq:wSO}
w=\tw_0\tx_0+\tw'_0\tx'_0+
\sum_{i=1}^{n-2}(w_ix_i+\tw_i\tx_i+w'_ix'_i+\tw'_i\tx'_i)
\end{equation}
Then   by writing out $w$  in matrix form and computing as in
(\ref{eq:Pcoeff}) and (\ref{eq:Ad})  we find the following formula for $P$: 
\begin{equation}\label{eq:PS0}
\begin{split}
P(w_i,\tw_i,w'_i,\tw'_i,\tw_0,\tw'_0)=&
{\textstyle \frac{1}{4}\left((
\sum_{i=1}^{n-2}w_iw_i'-\tw_i\tw_i')-\tw_o\tw'_0\right)^2}\\[8pt]
&{\textstyle+\left((\sum_{i=1}^{n-2}w_i'\tw_i)+\frac{1}{2}\tw_0^2\right)
\left((\sum_{i=1}^{n-2}w_i\tw_i')-\frac{1}{2}(\tw'_0)^2\right)}\\
\end{split}
\end{equation}    

Then plugging this expression for $P$ into (\ref{eq:fe}) we get
\begin{equation}
\label{eq:feSO}
f_{\psi}=\frac{1}{4f_0^3}(a^2+4bc-(f_0f_0')^2).
\end{equation} 
where
\begin{equation}
\label{eq:abcSO}
\begin{array}{l}
a=(\sum_{i=1}^{n-2}f_if_i'-\wt{f}_i\wt{f}_i')-\wt{f}_0\wt{f}_0'\\[12pt]
b=(\sum_{i=1}^{n-2}f_i\wt{f}_i')+\frac{1}{2}(\wt{f}_0)^2\\[12pt]
c=(\sum_{i=1}^{n-2}f_i'\wt{f}_i)-\frac{1}{2}(\wt{f}_0')^2
\end{array}
\end{equation}

By   making the substitutions, where $i=1,\dots,m$,  
\[f_i^\ep\mapsto\Th^{x_i^\ep},\qquad \tf_i^\ep\mapsto\Th^{\tx_i^\ep},
\qquad  \tf_0^\ep\mapsto\Th^{\tx_0^\ep}   \]
we get symbols $A,B,C$  corresponding to  the functions  $a,b,c$.
Then  plugging our expression for $P$ into  (\ref{eq:rzero}) we get
\begin{equation}
\label{eq:rzeroSO}
r_0=\frac{1}{4f_0}\left(A^2+4BC-\la^2(\Phi^{x_0})^2\right)
\end{equation}

The Weyl quantization of these principal symbols $A,B,C$ are the operators
\begin{equation}\label{ASO}
\begin{split}
\bA=&{\textstyle \frac{1}{2}\left(\sum_{i=1}^{n-2}
\Delta^{x_i}\Delta^{x'_i}+\Delta^{x'_i}\Delta^{x_i}-
\Delta^{\tx_i}\Delta^{\tx'_i}-\Delta^{\tx'_i}\Delta^{\tx_i}\right)}\\[8pt]
&{\textstyle -\half\left(\Delta^{\tx_0}\Delta^{\tx'_0}+\Delta^{\tx'_0}
\Delta^{\tx_0}\right)}\\
\end{split}
\end{equation}

\begin{equation}
\label{BSO}
{\bB}={\textstyle \left(\sum_{i=1}^{n-2}\Delta^{x_i}\Delta^{\tx'_i}\right)+
\half(\Delta^{\tx_0})^2}
\end{equation}
\begin{equation}
\label{CSO}
{\bC}={\textstyle \left(\sum_{i=1}^{n-2}\Delta^{\tx_i}\Delta^{x'_i}\right)-
\half(\Delta^{\tx'_0})^2}
\end{equation}

\begin{lem}
\label{lem:ABCSO}
$\bA$, $\bB$ and $\bC$  are differential operators on 
$\O\reg$ of order $2$ with respective principal symbols $A$, $B$, $C$.
$\bA$, $\bB$ and $\bC$ are all  Euler homogeneous of degree $0$ 
and each one commutes with the vector fields $\eta^z$, $z\in\heis$.  We have 
$[\eta^h,\bA]=-2\bA$, $[\eta^h,\bB]=-2\bB$, and $[\eta^h,\bC]=-2\bC$.
\end{lem}
{\bf Proof:} Same as for  Lemma \ref{lem:ASL}.\Qed

Now  routine calculation results in the formulas

\begin{equation}
\label{eq:ABCt_k}
{\bA}(f_{\psi}^k)=-t_k\frac{a}{f_0}f_{\psi}^{k-1}\qquad
{\bB}(f_{\psi}^k)=-t_k\frac{b}{f_0}f_{\psi}^{k-1}\qquad
{\bC}(f_{\psi}^k)=-t_k\frac{c}{f_0}f_{\psi}^{k-1}
\end{equation}
where

\begin{equation}
\label{eq:t_k}
t_k=k(k+\frac{m}{2}-1)
\end{equation}

Next more computation gives the following formulas, where we have set
$r_k=k+{\textstyle\frac{m}{2}}-2$.
 
\begin{equation}
\label{eq:A2SO}
{\bA}^2(f_{\psi}^k)=\left((k-1)\frac{a^2+4bc+r_ka^2}{f_0^2}+
mf_0\right)t_kf_{\psi}^{k-2}
\end{equation}

\begin{equation}
\label{eq:CBSO}
2{\bC\bB}(f_{\psi}^k)=\left((k-1)\frac{a^2+4bc+2r_kbc+af_0'f_0}{f_0^2}+
mf_0f_{\psi}\right)t_kf_{\psi}^{k-2}
\end{equation}

\begin{equation}
\label{eq:BCSO}
2{\bB\bC}(f_{\psi}^k)=\left((k-1)\frac{a^2+4bc+2r_kbc-af_0'f_0}{f_0^2}+
mf_0f_{\psi}\right)t_kf_{\psi}^{k-2}
\end{equation}

Now combining the last three calculations we get
\begin{equation}
\label{eq:combSO}
\begin{split}
&(\bA^2+2\bB\bC+2\bC\bB)(f_{\psi}^k)\\
&=\left[3mt_kf_0f_{\psi}+(k-1)(k+{\textstyle\frac{m}{2}}+1)t_k\frac{a^2+4bc}
{f_0^2}\right]f_{\psi}^{k-2}\\
\end{split} 
\end{equation}

Our aim is again to  find $q(E)$ so that
we can build the operator $\bS$  with the desired properties.
Proceeding as in the previous subsection,  we compare (\ref{eq:combSO}) with
(\ref{eq:etasq}). By (\ref{eq:feSO}) we have the relation
\begin{equation}
\label{eq:keySO}
\frac{a^2+4bc}{f_0^2}-(f_0')^2=4f_0f_{\psi}
\end{equation}
This is the unique linear relation in $R(\Oreg)$ between the functions
$\frac{a^2+4bc}{f_0^2}$, $(f_0')^2$, and $f_0f_{\psi}$.
As in the  previous subsection, we obtain
\begin{lem} 
\label{lem:mainSO}
Let 
\begin{equation}
\label{eq:SSO}
\bS=\frac{1}{4}\left((\bA^2+2\bB\bC+2\bC\bB)-q(E)(\eta^{x_0})^2\right)
\end{equation}
where $q(E)$ is some polynomial in the Euler operator $E$.   Then $\bS$ is a 
differential operator on $\O\reg$ of order $4$ with principal symbol $S$. 
$\bS$ is homogeneous  of degree $0$ and satisfies the conditions in
(\ref{eq:Sreq}).  The  polynomial  
\begin{equation}
\label{eq:qSO}
q(E)=(E+{\textstyle\frac{m}{2}}+1)(E+{\textstyle\frac{m}{2}}-1) 
\end{equation}
is the unique polynomial in  $E$ such that  $\bS$ satisfies (iii) in  Proposition 
\ref{prop:extend}. Then $\bS(f_{\psi}^{k})=\beta_kf_0f_{\psi}^{k-1}$
where
\begin{equation}
\label{eq:betaSO}
\beta_k=k(k+1)(k+\frac{m}{2}-1)(k+\frac{m}{2}-\half)
\end{equation}
\end{lem}
  
The polynomial $\beta_k$  has leading term $k^4$ and all coefficients 
non-negative  since $m\ge 2$. 
 Proposition \ref{prop:extend} now says that
$f_0^{-1}\bS$  is a  lowest weight vector of a copy of $\g$  in $\D_{-1}(\O)$ 
Thus we have proven

\begin{prop} 
\label{prop:SO} Let $\g=\so(m+4,\C)$, $m\ge 2$. The differential operator 
on $\Oreg$
\begin{equation}D_0=\frac{1}{f_0}\bS\end{equation}  
defined  by (\ref{eq:SSO}) and (\ref{eq:qSO})
extends to an algebraic differential operator  on $\O$.
$D_0$ satisfies  all three conditions in Theorem \ref{thm:quant}. We have 
\begin{equation}D_0(f_{\psi}^{k})=\beta_kf_{\psi}^{k-1}\end{equation} where 
$\beta_k$ is given by  (\ref{eq:betaSO}).
\end{prop}
 
\begin{rem}\rm
We can show that  $\bA$, $\bB$ and $\bC$  extend to  differential operators on
$\O$.
\end{rem}

\section{Relations to other work}
\label{sec:last}

We return to the more general situation of $K$-orbits on $\p$ discussed in
the introduction. The invariant theory of the action of $K$ on $\p$ was 
analyzed in \cite{kr}.

Suppose  $\g$ is a simple Lie algebra 
and  $\k$ has non-trivial center. Then the action of the center  of $K$ defines a 
$K$-invariant splitting  $\p=\p^+\oplus\p^-$. Let $Y$ be a $K$-orbit inside
$\p^+$. Then $\D(Y)$ has been studied for  various cases of this sort   
in \cite{lsm}, \cite{lsmst}, \cite{lst}  where the
functions and vector fields on $Y$ do not generate $\D(Y)$.   They construct 
an  extension of the natural infinitesimal action   $\k\to\Vect\, Y$ 
to a Lie algebra homomorphism $\pi:\g\to\D(Y)$ where $\p^+$ acts by
multiplication operators and 
$\p^-$ acts by  commuting order $2$ differential operators on $Y$
which are  Euler homogeneous of degree $-1$.

In particular,  in  \cite{lsmst},  Levasseur, Smith and Stafford    
analyzed cases  where   $Y=\Omin\cap\p^+$ and  $\Omin$ is
the minimal nilpotent orbit of  $\g$.
They  prove that the kernel of the algebra homomorphism $\U(\g)\to\D(Y)$  
defined by  their map $\pi$  is the  Joseph ideal.  In (\cite{lst}) 
Levasseur and Stafford analyzed $\D(Y)$  where  $Y$  is 
any  $K$-orbit in $\p^+$ and $\g$ is classical.  Using the theory of Howe 
pairs and the Oscillator representation, they  construct $\pi$ for $Y$
``sufficiently small" and   prove that the kernel of the corresponding algebra 
homomorphism $\U(\g)\to\D(Y)$ is a completely prime maximal ideal.

Now suppose $\g$ is simple and $\k$ has trivial center. Then $\p$ is 
irreducible as a $\k$-module. In \cite{bkDiff}, R. Brylinski and B. Kostant 
studied   cases of this sort where
$Y=\Omin\cap\p$. They constructed a space $L$ of   commuting
differential operators on $Y$ such that $L$ is $G$-irreducible and carries the 
$K$-representation $\p$. The operators in $L$ have order $4$ and  are Euler
homogeneous of degree $-1$. The highest weight vector in $L$ was constructed
as the quotient by a multiplication operator 
of a  homogeneous quartic polynomial in certain commuting
vector fields of the $\k$-action. 

Our methods are an extension of that  construction.  
However,   our quantization is considerably
more subtle as the our  principal symbol $r_0$ 
involves   a homogeneous quartic polynomial in the  symbols of {\it  non-commuting}
vector fields (which lie outside the $\k$-action).  The  quantization then requires not
only symmetrization but  also the introduction of  lower order correction terms. 
Theorem \ref{thm:quant} says  that these correction terms are uniquely 
determined.

In \cite{lo},  Lecomte and Ovsienko show that the algebra of polynomial 
symbols  on $\R^n$ admits a unique quantization (involving non-obvious lower
order terms) equivariant under the vector field action of
$\sl(n+1,\R)$ on $\R^n$ which arises by  embedding $\R^n$ as the  ``big
cell" in $\R{\mathbb P}^n$.   It would be interesting to see if there is a
$\g$-equivariant quantization map   ${\mathcal Q}:R(T^*\O)\to\D(O)$.


\vskip 4pc   


\begin{thebibliography}{10}


\bibitem[A-B1]{abr1}
{\bf A. Astashkevich, R. Brylinski} 
\newblock {\it Cotangent bundle models of complexified small  nilpotent  orbits}
\newblock preprint

\bibitem[A-B2]{abr2}
{\bf A. Astashkevich, R. Brylinski}
\newblock {\it  Lie algebras of  exotic pseudo-differential symbols}
\newblock in preparation

\bibitem[A-B3]{abr3}
{\bf A. Astashkevich, R. Brylinski}
\newblock {\it Quantization of complex minimal nilpotent orbits and construction of
algebraic star products}
\newblock in  preparation

\bibitem[Bo-Br]{bobr}{\bf W. Borho and J.-L.  Brylinski} {\it
Differential operators on homogeneous spaces I}, Invent. Math. 69 (1982),  
437-467



\bibitem[B1]{brInstantons}{\bf R. Brylinski} {\it
Instantons and Kaehler Geometry of Nilpotent Orbits},
in "Representation Theories and Algebraic Geometry", A. Broer ed., Kluwer,
85-125, in press.

\bibitem[B2]{brGQmin}{\bf R. Brylinski} {\it
Geometric quantization of real minimal nilpotent orbits}, Differential
Geometry and its Applications,  vol. 9/1-2  (1998),  5--58, in press.

\bibitem[B3]{brISC}{\bf R. Brylinski},   in  preparation.

\bibitem[B-K1]{bkHam}{\bf R. Brylinski, B. Kostant} {\it
Nilpotent orbits, normality and Hamiltonian group actions.} 
J. Amer. Math. Soc., vol.  7 (1994), no. 2, 269--298.

\bibitem[B-K2]{bkDiff}{\bf R. Brylinski, B. Kostant} {\it
Differential operators on conical Lagrangian manifolds.}
Lie theory and geometry, 65--96, Progr. Math., vol. 123, 
Birkhauser Boston, Boston, MA, 1994.

\bibitem[B-K3]{bkLagr}{\bf R. Brylinski, B. Kostant} {\it 
Lagrangian models of minimal representations of $E\sb 6$, 
$E\sb 7$ and $E\sb 8$.} Functional analysis on the eve of the 21st century, 
 Vol. 1 (New Brunswick, NJ, 1993), 13--63, Progr. Math., vol. 131, 
Birkh\"auser Boston, Boston, MA, 1995. 

\bibitem[Gar]{gar} {\bf D. Garfinkle}
\newblock {\it A new construction of the Joseph ideal},
\newblock MIT Doctoral Thesis, 1982

\bibitem[J]{jos} {\bf A. Joseph} {\it
The minimal orbit in a simple Lie algebra and associated maximal ideal.}
Ann. Scient. Ec. Norm. Sup., vol.  9 (1976), 1--30

\bibitem[K-R]{kr} {\bf  B. Kostant and S. Rallis} {\it
Orbits and representations associated with symmetric spaces},
Amer. J. Math, vol. 43 (1971), 753--809.

\bibitem[Kr]{kron} 
   {\bf P. B. Kronheimer}, {\it Instantons and the geometry of the nilpotent
   variety},  Jour. Diff. Geom.,   vol. 32 (1990), 473--490.
 


\bibitem[L-O]{lo} {\bf P.B.A. Lecomte and V. Ovsienko} {\it Projectively 
invariant symbol map and cohomology of vector fields Lie algebras intervening
in  quantization}, preprint


\bibitem[L-Sm]{lsm} {\bf T. Levasseur and S.P. Smith} {\it Primitive ideals and 
nilpotent orbits in type $G_2$}, J. Algebra, vol. 114 (1988), 81--105

\bibitem[L-Sm-St]{lsmst} {\bf T. Levasseur,  S.P. Smith and  J.T. Stafford} 
{\it The minimal nilpotent orbit, the Joseph ideal and differential operators},  
J. Algebra, vol. 116 (1988), 480--501
 

\bibitem[L-St]{lst} {\bf T. Levasseur and J.T. Stafford} {\it Rings of 
differential operators on classical rings of invariants}, Memoirs of the AMS, vol.
81, number 412, 1989

\bibitem[Se]{sek} {\bf   J. Sekiguchi} {\it Remarks on real nilpotent orbits of 
a symmetric  pair},  J. Math. Soc. Japan, vol. 39 (1987), 127--138.

\bibitem[Ve]{ver} {\bf  M. Vergne}, 
{\it Instantons et correspondance de Kostant-Sekiguchi}, 
 C.R. Acad. Sci. Paris, vol.    320 (1995), Serie 1, 901--906.


\bibitem[V-P]{pv} {\bf  E.B. Vinberg and V.L. Popov}, 
{\it  On a class of quasihomogeneous affine varieties}, 
 Math. USSR Izv., vol.  6  (1972 ),  743--758.
 

 
 
\end{thebibliography}
\end{document}